\documentclass[11pt]{amsart}
\usepackage{amsfonts,amsmath,amssymb}
\usepackage{url}
\usepackage{hyperref,amssymb}
\hypersetup{colorlinks=true,linkcolor=blue,
urlcolor=cyan,citecolor=magenta} 
\usepackage[hyperpageref]{backref}
\usepackage{frcursive,calrsfs,mathrsfs}
\usepackage{xcolor}
\newtheorem{theorem}{Theorem}[section]
\newtheorem{lemma}[theorem]{Lemma}
\newtheorem{corollary}[theorem]{Corollary}
\newtheorem{conjecture}[theorem]{Conjecture}
\newtheorem{proposition}[theorem]{Proposition}
\newtheorem{remark}[theorem]{Remark}
\newtheorem{remarks}[theorem]{Remarks}
\newtheorem{example}[theorem]{Example}
\newtheorem{definition}[theorem]{Definition}

\newcommand{\Q}{{\mathbb{Q}}}
\newcommand{\Z}{{\mathbb{Z}}}
\newcommand{\F}{{\mathbb{F}}}

\newcommand{\ds}{\displaystyle}

\newcommand{\ov}{\overline}
\newcommand{\wt}{\widetilde}

\newcommand{\ft}{\footnotesize}
\newcommand{\ns}{\normalsize}

\newcommand{\BH}{{\bf H}}
\newcommand{\CH}{{\mathcal H}}
\newcommand{\BE}{{\bf E}}
\newcommand{\CE}{{\mathcal E}}
\newcommand{\BF}{{\bf F}}
\newcommand{\CF}{{\mathcal F}}

\newcommand{\X}{{\bf X}}

\newcommand{\order}{\raise0.8pt \hbox{${\scriptstyle \#}$}}
\newcommand{\oorder}{\raise0.8pt \hbox{\tiny${\scriptstyle \#}$}}
\newcommand{\lien}{\mathrel{\mkern-4mu}}
\newcommand{\too}{\relbar\lien\rightarrow}
\newcommand{\tooo}{\relbar\lien\relbar\lien\too}

\newcommand{\plus}{\ds\mathop{\raise 0.5pt \hbox{$\bigoplus$}}\limits}
\newcommand{\prd}{\ds\mathop{\raise 1.0pt \hbox{$\prod$}}\limits}
\newcommand{\sm}{\ds\mathop{\raise 1.0pt \hbox{$\sum$}}\limits}
\newcommand{\ffrac}[2]{\hbox{\ft $\displaystyle\frac{#1}{#2}$}}
\newcommand{\limproj}{\mathop{\lim_{\longleftarrow}}}
\newcommand{\Gal}{{\rm Gal}}

\newcommand{\Norm}{{\bf N}}
\newcommand{\J}{{\bf J}}
\newcommand{\jj}{{\bf j}}
\newcommand{\Nu}{\hbox{\Large $\nu$}}
\newcommand{\Lbda}{{\bf \Lambda}}

\newcommand{\rk}{{\rm rank}}

\newcommand{\Ker}{{\rm Ker}}
\newcommand{\Coker}{{\rm Coker}}

\newcommand{\ram}{{\rm ram}}
\newcommand{\ab}{{\rm ab}}
\newcommand{\alg}{{\rm alg}}
\newcommand{\ar}{{\rm ar}}
\newcommand{\nr}{{\rm nr}}

\newcommand{\Hom}{{\rm H}}
\newcommand{\g}{\textbf{g}}
\def\plus{\displaystyle\mathop{\raise 2.0pt \hbox{$\bigoplus $}}\limits}
\def\prd{\displaystyle\mathop{\raise 2.0pt \hbox{$\prod$}}\limits}
\def\sm{\displaystyle\mathop{\raise 2.0pt \hbox{$\sum$}}\limits}

\author[Georges Gras]{Georges Gras}

\address{Villa la Gardette, chemin Ch\^ateau Gagni\`ere F-38520, 
Bourg d'Oisans  Url: {\rm \url{http://orcid.org/0000-0002-1318-4414}}}
\email{g.mn.gras@wanadoo.fr}



\begin{document}

\title[Chevalley--Herbrand formula and the Main Conjecture]
{The Chevalley--Herbrand formula \\ and the real abelian Main Conjecture \\
\ft {\bf New criterion using capitulation of the class group}\ns}

\date{September 2, 2022} 

\begin{abstract} 
The Main Theorem for abelian fields (often called Main Conjecture 
despite proofs in most cases) has a long history which has found a solution 
by means of ``elementary arithmetic'', as detailed in Washington's book 
from Thaine's method having led to Kolyvagin's Euler systems. 
Analytic theory of real abelian fields $K$ says (in the semi-simple case) 
that the order of the $p$-class group $\CH_K$ is equal to the $p$-index 
of cyclotomic units  $(\CE_K : \CF_K)$. We have conjectured (1977) the  
relations $\order \CH_\varphi = (\CE_\varphi : \CF_\varphi)$ for the isotypic 
$p$-adic components using the irreducible $p$-adic characters 
 $\varphi$ of $K$. We develop, in this article, new promising links between: 
(i) the Chevalley--Herbrand formula giving the number of ``ambiguous 
classes'' in $p$-extensions $L/K$, $L \subset K(\mu_\ell^{})$ for the 
auxiliary prime numbers $\ell \equiv 1 \pmod {2p^N}$ inert in $K$;
(ii) the phenomenon of capitulation of $\CH_K$ in $L$;
(iii) the real Main Conjecture $\order \CH_\varphi = (\CE_\varphi : \CF_\varphi)$ 
for all~$\varphi$. We prove that the real Main Conjecture is trivially fulfilled as soon 
as $\CH_K$ capitulates in $L$ (Theorem \ref{thmppl}). Computations 
with PARI programs support this new philosophy of the Main Conjecture. 
The very frequent phenomenon of capitulation suggests Conjecture \ref{conjcap}.
\end{abstract}

\maketitle

\tableofcontents

\section{Introduction -- Statement of the main result} 

\subsection{Abelian characters}
Let $\Q^\ab$ be the maximal abelian extension of $\Q$ contained in an 
algebraic closure $\ov \Q$ of $\Q$; let $\Q_p$ be the $p$-adic field and 
$\ov \Q_p$ an algebraic closure of $\Q_p$ containing $\ov \Q$. 

\smallskip
Let ${\bf \Psi}$ be the set of irreducible characters of $ \Gal(\Q^\ab/\Q)$, of degree~$1$ 
and finite order, with values in $\ov \Q_p$. We define the set ${\bf \Phi}$ of irreducible 
$p$-adic characters and the set $\X$ of irreducible rational characters.
For a subfield $K$ of $\Q^\ab$, we define the subsets 
${\bf \Psi}_K$, ${\bf \Phi}_K$, $\X_K$, whose kernels fixe $K$. The set $\X$ is in 
one-to-one correspondence with the set of cyclic subfields of $\Q^\ab$.

\smallskip
The notation $\psi \mid \varphi \mid \chi$ (for $\psi \in {\bf \Psi}$, $\varphi \in {\bf \Phi}$,
$\chi \in \X$) means that $\varphi$ is a term of $\chi$ and $\psi$ a term of $\varphi$;
so $\varphi$ is the sum of the $\Q_p$-conjugates of $\psi$ and $\chi$ the sum
of the $\Q$-conjugates of $\psi$ (cf. \cite{Ser1998}).

\smallskip
Let $\chi \in {\bf \X}$; we denote by $g_\psi = g_\varphi =
g_\chi$ the order of $\psi \mid \varphi \mid \chi$; the field of values of these
characters is the cyclotomic field $\Q(\mu_{g_\psi})$.

\smallskip
Let $K/\Q$ be an abelian extension of Galois group $g$ of prime-to-$p$ order.
For any $\Z[g]$-module $A$ of finite type, we denote by ${\mathcal A}_\chi$
(resp. ${\mathcal A}_\varphi$) the $\chi$-component (resp. the
$\varphi$-component) of ${\mathcal A} := A  \otimes \Z_p$; we get 
${\mathcal A}_\chi = \bigoplus_{\varphi \mid \chi}{\mathcal A}_\varphi$,
then ${\mathcal A} = \bigoplus_{\chi \in \X_K} \bigoplus_{\varphi \mid \chi}{\mathcal A}_\varphi
= \bigoplus_{\varphi \in {\bf \Phi}_K}{\mathcal A}_\varphi$.

\smallskip
In this semi-simple case all reasonings reduce to the arithmetic of the 
cyclic subfields of $K$; so in what follows, we only consider cyclic $K/\Q$'s. 

\subsection{Main theorem}
We have obtained the following result (Theorem \ref{relfond}):

\begin{theorem} \label{thmppl}
Let $p \geq 2$ be a prime. Let $K/\Q$ be a real cyclic extension of prime-to-$p$ degree 
and Galois group $g$. Denote by $\CH_K = \bigoplus_{\varphi \in {\bf \Phi}_K} \CH_\varphi$ 
the $p$-class group of $K$.

\smallskip\noindent
Consider primes $\ell \equiv 1 \pmod {2 p^N}$, $\ell$ totally inert in $K$, and 
let $K_n$ be the subfield of $K(\mu_\ell^{})$ of degree $p^n$ over $K$, 
$n \in [1,N]$, where $\mu_\ell^{}$ is the group of $\ell$-roots of unity.

\smallskip\noindent
Let $\BE_K$ (resp. $\BF_K$) be the group of units (resp. of 
cyclotomic units) of $K$ and put $\CE_K = \bigoplus_{\varphi \in {\bf \Phi}_K} \CE_\varphi$ 
(resp $\CF_K = \bigoplus_{\varphi \in {\bf \Phi}_K} \CF_\varphi$). 

\smallskip
Then:

\smallskip
(i) As soon as $\CH_K$ capitulates by extension in some $K_n$, the Main Conjecture 
holds in $K$, that is to say, $\order \CH_\varphi = (\CE_\varphi : \CF_\varphi)$
for all $\varphi \in {\bf \Phi}_K$.

\smallskip
(ii) We have $\order \CH_{K_n} = \order \CH_K$, for all $n \in [1, N]$, if and only if
$\order \CH_{K_1} = \order \CH_K$. 
If this stability property holds, we have the following consequences:

\smallskip
\qquad $\bullet$ The $p$-class groups $\CH_{K_n}$ are invariant by 
$\Gal(K_n/K)$ and the norms $\Norm_{K_n/K} : \CH_{K_n} \too \CH_K$ are isomorphisms.

\smallskip
\qquad $\bullet$ Let $p^e$ be the exponent of $\CH_K$ and assume 
$N \geq e$; then the $p$-class group $\CH_K$ capitulates in $K_e$ and (from (i))
the Main Conjecture holds.
\end{theorem}

Denote by $\Norm$ the arithmetic norm and by $\J$ the transfer 
(corresponding to extension of classes).
Note, once for all, that if $\CH_K$ is of exponent $p^e$ and capitulates in $K_n$, 
the relation $\Norm_{K_n/K} \circ \J_{K_n/K}(\CH_K) = \CH_K^{p^n}$ shows than 
necessarily $n \geq e$; but incomplete capitulation may occur for any $n \geq 1$.

\smallskip
We have proposed in \cite{Gra2022$^c$} the following Conjecture of capitulation, 
without any assumption of splitting on $\ell$:

\begin{conjecture}\label{conjcap}
Let $K$ be any totally real number field and let $\CH_K$ be its $p$-class group,
of exponent $p^e$. 
There exist infinitely many primes $\ell \equiv 1 \pmod {2 \,p^N}$, $N\geq e$, such that 
$\CH_K$ capitulates in $K(\mu_\ell^{})$.
\end{conjecture}

We shall give Section \ref{comput} extensive numerical computations with
PARI programs \cite{Pari2013} showing that capitulation in such auxiliary cyclic
$p$-extensions is, surprisingly, very frequent and conjecturally holds for infinitely many $\ell$. 
The criterion (ii) of the theorem allows easy effective verifications; but capitulation
may hold without stability in $K_1/K$.

\smallskip
We do not intend to evoke the case of the abelian capitulations of class 
groups proved in the literature (Gras \cite{Gra1997}, Kurihara \cite{Kur1999}, 
Bosca \cite{Bos2009}, Jaulent \cite{Jau2022}); all techniques in these papers need 
to built abelian $p$-extensions $L_0$ of $\Q$, ramified at various primes 
and requiring many local arithmetic conditions, whose compositum $L$ with $K$ gives 
a capitulation field of $\CH_K$; the method is completely incomparable to ours since 
it must apply to any real abelian field $K$, of arbitrary degree, obtained in an iterative 
process giving that the maximal real subfield of $\Q\big( \bigcup_{f>0} \mu_f^{} \big)$ 
is principal.

\smallskip
However, these results together with Theorem \ref{thmppl} and the help of 
numerical computations, among many other results of the literature, support the 
fact that the phenomenon of capitulation governs many aspects of abelian arithmetic, 
independently of the well-known case of capitulation in the Hilbert class field, from 
Hilbert's theorem 94 and a lot of improvements (see the surveys 
\cite{Jau1988,Jau1998,Mai1997,Mai1998} and their references).

\smallskip
In a slightly different, but related, context of the $p$-adic class field theory, mention for 
instance that the capitulation of the logarithmic class group of $K$ \cite{Jau1994,Jau1998}, in its 
cyclotomic $\Z_p$-extension $K_\infty = \bigcup_{n \geq 0} K_n$, is equivalent to 
Greenberg's conjecture \cite{Gree1976} saying that the Iwasawa invariants $\lambda$, $\mu$
for $\ds \limproj \CH_{K_n}$ are zero \cite{Jau2016,Jau2019$^a$,Jau2019$^b$}. 

\smallskip
Analogous criteria of stability of the $\order \CH_{K_n}$ were given by Greenberg  
then by Fukuda \cite{Fuk1994}; for instance, in the similar context as ours where one 
considers $p$ non split in $K$, the condition $\lambda = \mu =0$ in $L=K_\infty$ is 
equivalent to the capitulation of $\CH_K$ in some $K_{n_0}$ \cite[Theorem 1]{Gree1976};
the case of $p$ totally split in $K$ \cite[Theorem 2]{Gree1976} relies on equalities
$\CH_{K_n}^{G_n} = \CH_{K_n}^\ram$ involving the exact sequence of Theorem \ref{suitecap}.

\smallskip
In a numerical setting, let's point out works of Kraft--Schoof \cite{KrSc1995$^a$, KrSc1995$^b$} 
(resp. Pagani \cite{Pag2022}), verifying Greenberg's conjecture for some real quadratic fields of 
conductor $f<10^4$ and $p=3$ (resp. $p=2$), by means of the analytic formulas in
some layers of $K_\infty$ with the use of cyclotomic units, giving some algebraic similarities.

\smallskip
But Greenberg's conjecture is still unproved and depends (without any assumption 
on the decomposition of $p$) on random algorithmic process, as explain in \cite{Gra2021}, 
governed by the torsion group of the Galois group of the maximal abelian $p$-ramified 
pro-$p$-extension of $K$ (essentially, the second Tate--Chafarevich group of $K$); 
this takes place in a deep $p$-adic context, beyond Leopoldt's conjecture.

\begin{remark}\label{stab}
{\rm We have proven the criterion (ii) of stability of Theorem \ref{thmppl} in 
\cite[Theorem 3.1, Corollaire 3.2]{Gra2022$^c$}, generalizing similar
results \cite{Fuk1994,LOXZ2022,MiYa2021}. 
More precisely, this criterion can be applied at some layer $n_0$ and one obtains
$\order \CH_{K_n} = \order \CH_{K_{n_0}}$ for all $n \geq n_0$, if and 
only if the equality holds for $n=n_0+1$; so, this means that $\CH_{K_{n_0}}$
capitulates in $K_{n_0+e_{n_0}}$, where $p^{e_{n_0}}$ is the exponent
of $\CH_{K_{n_0}}$, but a fortiori, $\CH_K$ capitulates in $K_{n_0+e_{n_0}}$.
We shall give again a proof in our particular simpler framework 
(Theorem \ref{criterion}).}
\end{remark}

\subsection{Methodology}
We shall obtain Theorem \ref{thmppl} by means of a classical exact sequence
describing $\CH_L^G$, in cyclic $p$-extensions $L/K$ of Galois group $G$, 
in terms of the units and image of the extension $\J_{L/K} : \CH_K \to \CH_L^G$ of $p$-classes, 
which gives, under the phenomenon of capitulation, the needed information about the 
index of cyclotomic units, taking into account their norm properties in abelian 
extensions (Corollary \ref{maincoro} and Proposition \ref{cycloformula}).
For that, we shall need the order of the ``$\varphi$-components'', for all $\varphi \in \Phi_K$, 
of the Chevalley--Herbrand formula giving $\order \CH_L^G$, that is to say, the 
computation of $\order \CH_\varphi^G$, where $\CH_\varphi^G := (\CH_\varphi)^G =
(\CH^G)_\varphi$.
This will give the opportunity to write cohomological exact sequences linking 
invariant classes and capitulation to the norm properties of the units (local and global);
for this, we shall follow \cite[III, p. 167]{Jau1986}. This step is crucial
since Chevalley--Herbrand formula writes, for $p$-class groups in $L/K$ cyclic
real of degree $p^n$, as follows \cite[pp. 402-406]{Che1933}:
\[\order \CH_L^G =  \ffrac{\order \CH_K \times \prod_{\mathfrak q} e_{\mathfrak q}(L/K)}
{p^n \times (\CE_K : \CE_K \cap \Norm_{L/K}(L^\times))}, \] 
where $e_{\mathfrak q}(L/K)$ is the ramification index in $L/K$ of the prime ideals 
${\mathfrak q}$ of $K$; the original Chevalley formula depends on the Herbrand 
theorem \cite{Her1930}, introducing the general ``Herbrand quotient'' (see Lemme 3, 
p. 375 of Chevalley's Thesis, then \cite[Appendice, \S\,1, p. 57]{Her1936}).
\,\footnote{It may be useful to recall some historical context of such 
an important formula, quoting the following excerpt from Chevalley's Thesis 
(footnote 19, p. 402): 
{\it Le calcul fait ici est la g\'en\'eralisation au cas ``cyclique quelconque'' du calcul 
fait par Takagi dans le cas ``cyclique de degr\'e premier'', calcul dont l'id\'ee essentielle 
se trouve d\'ej\`a dans le ``Zahlbericht'' de Hilbert, dans la d\'emonstration du th\'eor\`eme 
suivant~: si un sur-corps relativement cyclique de degr\'e premier de $k$ est non ramifi\'e, 
il y a au moins un id\'eal de $k$ qui n'est pas principal dans $k$ mais qui est principal 
dans le sur-corps. L'extension au cas ``cyclique de degr\'e quelconque'' a \'et\'e rendue 
possible par le th\'eor\`eme des unit\'es de Herbrand}.}

\smallskip
One must obtain the decompositions of each factor into $\varphi$-components, especially 
for $\order \CH_K$ since $\CH_K$ {\it is not} a submodule of  $\CH_L^G$, the transfer 
map $\J_{L/K} : \CH_K \to \CH_L^G$ being in general non injective (see, e.g., 
Example \ref{example1}). For this, we introduce the following elementary principle.

\begin{definition} 
Let $K/\Q$ be an abelian extension of Galois group $g$
and let $X$ be a finite $\Z[g]$-module. Assume that we know that
$\order X = \prod_i \order U_i \times \prod_j (\order V_j)^{-1}$ depending on
finite $\Z[g]$-modules $U_i,V_j$. 

\smallskip\noindent
We say that $X$ (or the formula giving $\order X$) is $p$-localizable 
if there exist exact sequences $1\to A_k \to B_k \to C_k \to 1$ of finite
$\Z[g]$-modules of the form $U_i$, $V_j$ and $X$,
such that $\order X$ is obtained by means
of the relations $\order B_k = \order A_k \times \order C_k$.
\end{definition}

Under this property, the flatness of $\Z_p$ allows to deduce the family of exact 
sequences of $\Z_p[g]$-modules $1\to {\mathcal A}_k := A_k \otimes \Z_p \to 
{\mathcal B}_k := B_k \otimes \Z_p  \to {\mathcal C}_k := C_k \otimes \Z_p  \to 1$, then,
taking the isotopic components (for instance by means of the fundamental
idempotents $e_\varphi$ of $\Z_p[g]$, $\varphi \in \Phi_K$), we get the exact sequences
$1\to {\mathcal A}_{k,\varphi} \to {\mathcal B}_{k,\varphi}\to 
{\mathcal C}_{k,\varphi} \to 1$
and the formulas $\order {\mathcal B}_{k,\varphi} = \order {\mathcal A}_{k,\varphi} 
\times \order {\mathcal C}_{k,\varphi}$ yielding $\order X_\varphi$. 

\smallskip
For instance, cohomology groups $\Hom^n (G,X)$, $n \in \{1,2\}$, $G =:\langle \sigma \rangle$,
are $p$-localizable, as soon as $X$ is $p$-locali\-zable; indeed:
$$\hbox{$\Hom^1(G,X) := \Ker_X(\Nu_{\!G})/X^{1-\sigma}$\  and \  
$\Hom^2(G,X) := X^G/\Nu_{\!G}(X)$,}$$
where $\Nu_{\!G}$ is the algebraic norm.

\subsection{The real abelian Main Conjecture}

The Main Conjecture for real abelian fields $K$ (to be called Main Theorem because of its 
numerous proofs; but this name has become common) essentially says that, when 
$p \nmid \order g$, one has, for  all $p$-adic irreducible characters $\varphi \in {\bf \Phi}_K$
\cite{Gra1976,Gra1977}:
\[\order \CH_{K,\varphi} = (\CE_{K,\varphi} : \CF_{K,\varphi}), \] 
for the $\varphi$-components of $\CH_K$, where $\CE_K = 
\BE_K \otimes \Z_p$, $\CF_K = \BF_K \otimes \Z_p$ denote the 
groups of global units and of cyclotomic units of $K$, respectively.

\smallskip
The  following obvious property of rational characters is to be considered as 
the ``Main Theorem'' for rational components \cite[Chap. I, \S\,1, 1]{Leo1954}:

\begin{theorem} \label{chiformula}
Let $K/\Q$ be an abelian extension; let $(A_\chi)_{\chi \in \X_K}$,
$(A'_\chi)_{\chi \in \X_K}$, be two families of positive numbers, indexed by the 
set $\X_K$ of irreducible rational characters of $K$. If for all subfields $k$ of $K$, 
one has $\prod_{\chi \in \X_k} A'_\chi = \prod_{\chi \in \X_k} A_\chi$,
then $A'_\chi = A_\chi$ for all $\chi \in \X_K$.
\end{theorem}

So this result applies to the well-known complex analytic global formula
$\order \BH_K = (\BE_K : \BF_K)$ for $K$ cyclic real, which implies  
$\order \CH_{K,\chi} = (\CE_{K,\chi} : \CF_{K,\chi})$
for all $\chi \in \X_K$, in the semi-simple case. But when $\chi$ is a sum of several 
$p$-adic irreducible characters $\varphi$, the complex analytic theory
does not give precise relation between $\order \CH_{K,\varphi}$ and 
$(\CE_{K,\varphi} : \CF_{K,\varphi})$.
 In other words, there is no sufficient information, even with the help of
 $p$-adic $\zeta$ and $L$-functions.

\smallskip
Deep geometrical methods (Ribet, Mazur--Wiles) were successful, then more 
arithmetic ones were used to solve the problem for odd and even characters
\cite{CoSu2006, Grei1992, Kol2007, PeRi1990, Rib2008, Rub1990, Thai1988} 
among others. The beginnings of the story around the proof, by Ribet \cite{Rib1976}, 
of the (non-trivial) converse of the Herbrand theorem \cite{Herb1932} on Bernoulli's 
numbers and the minus part of the $p$-class group of the $p$th cyclotomic field, 
is given in the survey \cite{Rib2008}.

\smallskip
A more complete and recent story is available in Washington's book 
\cite[Chap.\,8]{Was1997} in which Thaine's paper is reproduced before 
the presentation of the various developments.

\smallskip
The non-semi-simple case of the real abelian Main Conjecture does exist and 
may be be split into two frameworks:

\smallskip
(i) The Iwasawa formulation \cite{Iwa1964}, ``replacing'' $K$ by its cyclotomic 
$\Z_p$-extension $K_\infty$ and considering semi-simple Galois actions from 
base fields of prime-to-$p$ degree, also called Main Conjecture without any precision 
(Lang--Rubin \cite{Lan1990}, Greither \cite{Grei1992}, Washington's book 
\cite[\S\S\,13.6, 15.4]{Was1997}); a typical statement being, for $K = \Q(\mu_p^{})$
and $\ds X := \limproj \CH_{K_n}$, the equality ${\rm char}(X_\varphi) = 
f(T,\varphi^*)\cdot u(T)$, with $\varphi \in \Phi_K$ even, $u(T)$ invertible in 
$\Z_p[[T]]$ and the power series $f(T,\varphi^*)$ such that 
$L_p(s,\varphi^*) = f((1+p)^s-1,\varphi^*)$ in terms of $p$-adic $L$-functions,
where $\varphi \mapsto \varphi^*$ by reflection.
Nevertheless, in the real case, 
Greenberg's conjecture makes it somewhat unnecessary and brings back to 
the finite cases $K_n/K$ as we have explained.

\smallskip
(ii) The case of cyclic extensions $K/\Q$ when $p \mid [K : \Q]$ and $K = K_\chi$; 
this case corresponds to our conjecture given in \cite{Gra1976,Gra1977} and still 
unproved for real fields.
We refer to the survey \cite{Gra2022$^a$} devoted to this non semi-simple case 
using the specific notion of $\varphi$-objects that we had introduced in the 1976's.
Indeed, classical works deal with an algebraic definition of the
$\varphi$-components of $p$-class groups, denoted $\CH^\alg_{K,\varphi}$,
which presents an inconsistency regarding analytic formulas; 
that is to say, when $g := \Gal(K/\Q)$ is cyclic of order $g_\chi \equiv 0 \pmod p$:
\[\CH^\alg_{K,\varphi} := \CH_K \otimes^{}_{\Z_p[g]} \Z_p[\mu_{g_\chi}],\ \,
\hbox{for all $\varphi \mid \chi$}, \]
with the $\Z_p[\mu_{g\chi}]$-action $\tau \in g \mapsto \psi (\tau)$
($\psi \mid \varphi$ of order $g_\chi$). We then have:
\[\CH^\alg_{K,\varphi} = \{x \in \CH_K, \ \Nu_{K/k}(x) = 1,\, 
\forall \, k \varsubsetneqq K \} \otimes^{}_{\Z_p[g]} \Z_p[\mu_{g_\chi}] \]
(where $\Nu_{K/k}$ is the algebraic norm), contrary to our arithmetic definition:
\[\CH^\ar_{K,\varphi} := \{x \in \CH_K, \ \Norm_{K/k}(x) = 1,\, \forall \, 
k \varsubsetneqq K \} \otimes^{}_{\Z_p[g]} \Z_p[\mu_{g_\chi}]. \]
This notion gives rise to an unexpected semi-simplicity
especially in accordance with analytic formulas, which enforce the conjecture in that~case:
\[\CH^\ar_{K,\chi}  := \{x \in \CH_K,\ \,\Norm_{K/k}(x) = 1,\, \forall \,
k \varsubsetneqq K \}  \otimes \Z_p = \plus_{\varphi \mid \chi} \CH^\ar_{K,\varphi} .\]

\section{The \texorpdfstring{$p$}{Lg}-localization of the Chevalley--Herbrand formulas}

\subsection{The higher rank Chevalley--Herbrand formulas}

Let $K/\Q$ be any extension. To avoid technical complications, we assume that
the given extension $L/K$, cyclic of degree $p^n$, $n \geq 1$, of Galois group 
$G =: \langle \sigma \rangle$, is such that any prime ideal of $K$, ramified in $L/K$, 
is totally ramified and that infinite places do not ramify (i.e., do not complexify) when $p=2$. 
Moreover we assume that the set of ramified primes is non-empty.

\smallskip
The basic principle of Chevalley--Herbrand formula for $p$-class groups
defines a filtration of the form $\{\CH_L^i\}_{i \geq 0}$, where $\CH_L^0 :=1$, 
$$\CH_L^1 := \CH_L^G, \ \hbox{of order} \   \order \CH_K \times \frac{p^{n (r-1)}}{(\CE_K : \CE_K 
\cap \Norm_{L/K}(L^\times))} ,$$ 
$\CH_L^{i+1}/\CH_L^i := (\CH_L/\CH_L^i)^G$,
up to a bound $m$ for which $\CH_L^m = \CH_L$ \cite{Gra2017}:
\begin{equation}\label{filtration}
\order (\CH_L^{i+1}/\CH_L^i) =  \frac{\order \CH_K}{\order \Norm_{L/K}(\CH_L^i)}
\!\times \!\frac{p^{n (r-1)}}{(\Lbda_K^i : \Lbda_K^i \cap \Norm_{L/K}(L^\times))},\  i \geq 0,
\end{equation}
where $r \geq1$ is the number of primes of $K$ ramified (totally) in $L/K$, and
$\Lbda_K^i := \{x \in K^\times,\ \, (x) = \Norm_{L/K}({\mathfrak A}),\, 
{\mathfrak A}P_L \in \CH_L^i \} \otimes \Z_p$ is a subgroup of finite type of $K^\times$ containing 
$\BE_K$ (with $\Lbda_K^0 = \BE_K$); the quotient 
$\Lbda_K^i /\Lbda_K^i \cap \Norm_{L/K}(L^\times)$ is of course a $p$-group
of order a divisor of $p^{n (r-1)}$. The first factor is called the class factor and the 
second one the norm factor.

\smallskip
In the case $r=1$ (which will be our context, taking a prime $\ell$ inert
in the cyclic extension $K/\Q$ and $L \subset K(\mu_\ell^{})$), one obtains only the class factors:
\begin{equation}\label{filtrationbis}
\order \CH_L^G = \order \CH_K  \ \ \ \& \ \ \ 
\order (\CH_L^{i+1}/\CH_L^i) =  \frac{\order \CH_K}{\order \Norm_{L/K}(\CH_L^i)},\  i \geq 1.
\end{equation}

There are a lot of generalizations of the original Chevalley--Herbrand formula. 
See for instance \cite[III, p. 167]{Jau1986}, \cite{Gra2017} with modulus,  
\cite{LiYu2020} in the spirit of Chevalley's theory of id\`eles, then
\cite{Gon2008} showing universality of the Chevalley--Herbrand
principle in arithmetic and algebraic geometry.

\subsection{The \texorpdfstring{$p$}{Lg}-localizations for \texorpdfstring{$\CH_L^G$}{Lg}}

We choose to privilege the aspect ``fixed points formula'', instead of genus theory 
(that is also $p$-localizable and may lead to parts of the results), to deduce 
(in \S\,\ref{ssect}) the $p$-localization of the filtration $(\CH_L^i)_{i \geq 1}$, 
which is of fixed points type, and to prove the stability Theorem \ref{criterion} 
using Chevalley--Herbrand formula in the $p$-tower $L/K$. We think that this 
is more suitable in a logical point of view since proofs of the deep statements 
of class field theory in the classical way rely first on Chevalley--Herbrand formula, 
whence Herbrand theorem, genus theory dealing only 
with local norms in non necessarily cyclic $p$-extensions 
for which Chevalley--Herbrand formula is not valid refering to {\it global norms}.
So,  to justify the integrality of some $p$-localized expressions,
one may use, at the end, the deep Hasse norm theorem saying that, in the cyclic 
case, $x \in K^\times$ is in $\Norm_{L/K}(L^\times)$ if and only if $x$ is everywhere 
local norm (except one place). All these comments about class field theory may 
be found in all the literature, and the classical way is detailed in our book
\cite[\S\,II.6; IV\,(b); Theorem II.6.2]{Gra2005}, respectively.

\subsubsection{Exact sequence of the ambiguous classes}
The $p$-localization of the Chevalley--Herbrand formula will exist from the definitions 
of class groups, ideal groups and units. Such $p$-localizations were given many 
years ago (\cite[Th\'eor\`emes II.1, II.2]{Gra1978}, \cite[Th\'eor\`emes III.1.12, 
III.1.13]{Jau1986}); these papers being written in french, we give again, in a more 
direct manner, the computations for the convenience of the reader.

\smallskip
To get a formula for the orders of the $\varphi$-components of $\CH_L^G$, 
we follow the process given in Jaulent's Thesis \cite[Chapitre III, p. 167]{Jau1986} 
(note that large generalizations of such methods, with ramification and decomposition, 
are also given in \cite{Mai1997, Mai1998} and, in the Galois case $L/K$, \cite{Gon2006} 
with many references):

\begin{theorem}
Let $L_0/\Q$ be a real cyclic extension of degree $p^n$, $n \geq 1$.
Let $K/\Q$ be a real abelian extension of Galois group $g$ and of
prime-to-$p$ degree. Put $L = L_0 K$ and $G := \Gal(L/K)
=: \langle \sigma \rangle$. We identify $\Gal(L/L_0)$ and $g$.

\smallskip\noindent
Let $I_K$ and $I_L$ (resp.  $P_K$ and $P_L$) be the ideal groups
(resp. the subgroups of principal ideals), of $K$ and $L$, respectively. 
Then, put $\BH_K := I_K/P_K$ and $\BH_L := I_L/P_L$.

\smallskip
(i) We have the exact sequence of $\Z[g]$-modules:

\smallskip
$1 \to \Ker(\J_{L/K}) \too  \Hom^1(G,\BE_L) \too \Coker (\jj_{L/K}^{}) \too $

\smallskip 
$\hspace{4cm}\Coker (\J_{L/K}) \too \Hom^2(G, \BE_L) \too  \Hom^2(G, L^\times)$,

\smallskip\noindent
where $\J_{L/K}$ is the transfer map $\BH_K \to \BH_L^G$ defined by
${\mathfrak a}P_K \mapsto {\mathfrak a}P_L$ and  $\jj_{L/K}^{}$ is
the extension of ideals $I_K \to I_L^G$.

\smallskip
(ii) Let $\varphi$ be an irreducible $p$-adic character of $K$. Then:

\smallskip
\qquad $\bullet$ If $\varphi = 1$, then $\order \CH_{L,\varphi}^G = 
\order \CH_{K,\varphi} = 1$;

\qquad $\bullet$ If $\varphi \ne 1$, then $\ds \order \CH_{L,\varphi}^G = 
\order \CH_{K,\varphi} \times  \frac{\order \Coker (\jj_{L/K}^{})_\varphi}
{(\CE_{K,\varphi} : \CE_{K,\varphi} \cap \Norm_{L/K}(L^\times))}$.
\end{theorem}

\noindent{\bf Proof.}
Note that each prime $\ell$, ramified in $L/K$, is totally ramified. 

\smallskip
Consider the exact sequences of $\Z[g]$-modules:
\begin{equation}\label{ab}
(a)\  1 \to \BE_L \to L^\times \to P_L \to 1, 
\ \ \ (b) \  1 \to \BE_K \to K^\times \to P_K \to 1
\end{equation}
\begin{equation}\label{cd}
(c)\  1 \to P_L \to I_L \to \BH_L \to 1,
\ \ \ (d) \ 1 \to P_K \to I_K \to \BH_K \to 1.
\end{equation}

\begin{lemma}\label{lemma0}
We have the following properties:

\smallskip
(i) $P_L^G/\jj_{L/K}^{}(P_K ) \simeq \Hom^1(G,\BE_L)$;

\smallskip
(ii) $\Hom^1(G,P_L) \simeq \Ker \big[  \Hom^2(G,\BE_L) \to \Hom^2(G,L^\times)  \big]
= \BE_K/ \BE_K \cap \Norm_{L/K}(L^\times)$;

\smallskip
(iii)  $\Hom^1(G,I_L) = 1$.
\end{lemma}

\noindent{\bf Proof.}
We have, from the above exact sequences \eqref{ab}\,(a),\,(b):
$$1 \to \BE_L^G = \BE_K \to L^\times{}^G = K^\times \to P_L^G  \to 
\Hom^1(G,\BE_L) \to  \Hom^1(G,L^\times) $$
$$\hspace{1.5cm} \to  \Hom^1(G,P_L) \to \Hom^2(G,\BE_L) \to \Hom^2(G,L^\times).$$

Since $\Hom^1(G,L^\times)  = 1$ (Hilbert's Theorem 90), this yields (i) and (ii).

\smallskip
The claim (iii) is classical since $I_L$ is a $\Z[G]$-module generated
by the prime ideals of $L$ on which the Galois action is canonical.
\qed

\medskip
From the above exact sequences \eqref{cd}\,(c),\,(d), we have the following 
commutative diagram:
$$\begin{array}{ccccccccc}
1  &  \!\! \tooo \!\! & P_K  &  \!\! \tooo \!\! & I_K &  \!\! \tooo \!\! & \BH_K &   \!\! \tooo \ \ 1 & \\   
\vspace{-0.3cm}   \\
& &   \Big \downarrow &  
&&\hspace{-1.3cm}\Big \downarrow \hbox{\ft$\jj_{L/K}^{}$\ns}  
&& \hspace{-2.0cm}\Big \downarrow  \hbox{\ft$\J_{L/K}$\ns}  \\   
 \vspace{-0.3cm}    \\
1  &   \!\! \tooo \!\!   & P_L^G  &  \!\! \tooo \!\! &  I_L^G &  \!\! \tooo \!\! &  \BH_L^G & 
 \!\! \!\! \!\! \!\! \!\! \tooo \!\! & \!\! \!\!  \!\! \!\! \!\Hom^1(G,P_L) \,.
\end{array} $$

The snake lemma gives the exact sequence:
$$\hspace{-2.8cm}1 \to  1  \too \Ker(\jj_{L/K}^{})  \too \Ker(\J_{L/K})  \too \Hom^1(G,\BE_L)   \too  $$
$$\hspace{2.0cm}\Coker(\jj_{L/K}^{}) \too \Coker(\J_{L/K})  \too  \Hom^2(G,\BE_L) 
\too \Hom^2(G,L^\times) ;$$
since $\Ker(\jj_{L/K}^{}) = 1$ and
$\Hom^2(G,L^\times) \simeq K^\times/\Norm_{L/K}(L^\times)$, it becomes:
$$\hspace{-2.8cm}1  \to  \Ker(\J_{L/K}) \to \Hom^1(G,\BE_L)   \to \Coker(\jj_{L/K}^{})  \to $$
$$\hspace{2.0cm}\Coker(\J_{L/K})  \to \Hom^2(G,\BE_L) 
\to \BE_K/\BE_K \cap \Norm_{L/K}(L^\times) \to 1.$$

\begin{remark}\label{coker0}
{\rm The representation $\Coker (\jj_{L/K}^{})$, as $g$-module, depends on the splitting 
of the prime ideals of $K$ ramified in $L/K$ and gives standard $\varphi$-components
from the relation $\order \Coker (\jj_{L/K}^{}) = \prod_{\mathfrak q} e_{\mathfrak q}(L/K)$; 
we refer to \cite{Gra1978} or \cite{Jau1986} for the most general localized formulas. 
From Hasse's norm theorem, the factors
$\ffrac{\order \Coker(\jj_{L/K}^{})_\varphi }
{(\CE_{K,\varphi} : \CE_{K,\varphi} \cap \Norm_{L/K}(L^\times))}$
are always integers in the totally ramified case.
In our particular context, where ramified primes $\ell$ are inert in $K/\Q$
and totally ramified in $L/K$, $\Coker (\jj_{L/K}^{})$ is of character $1$, 
so all the factors $\ffrac{\order \Coker(\jj_{L/K}^{})_\varphi }
{(\CE_{K,\varphi} : \CE_{K,\varphi} \cap \Norm_{L/K}(L^\times))}$ are trivial 
for $\varphi \ne 1$; in particular, any unit of $K$ is a global norm in $L/K$. 
Note that if the $p$-Hilbert class field $H_K^\nr$ of $K$ is not disjoint from $L/K$,
the Chevalley--Herbrand formula becomes $\order \CH_L^G = [H_K^\nr : L \cap H_K^\nr] 
\times \ffrac{\prod_{\mathfrak q} e_{\mathfrak q}(L/K)}{[L : L \cap H_K^\nr] \times 
(\CE_K : \CE_K \cap \Norm_{L/K}(L^\times))}$ giving another localization
(e.g., $L=H_K^\nr$ cyclic gives $\CH_L^G =1$).}
\end{remark}

We obtain the $p$-localized exact sequences of $\Z_p[g]$-modules:
$$\hspace{-2.8cm}1 \to  \Ker(\J_{L/K})_\varphi \to  \Hom^1(G,\CE_{L,\varphi}) \to  
\Coker(\jj_{L/K}^{})_\varphi  \to  $$
$$\hspace{1.5cm}\Coker(\J_{L/K})_\varphi  \to  \Hom^2(G,\CE_{L,\varphi}) \to 
\CE_{K,\varphi}/\CE_{K,\varphi} \cap \Norm_{L/K}(L^\times) \to 1$$
(by $p$-localization, $\Hom^1(G,\BE_L)_\varphi \!= \Hom^1(G,\CE_{L,\varphi})$,
$\Hom^2(G,\BE_L)_\varphi\! = \Hom^2(G,\CE_{L,\varphi})$).

\smallskip
All these $\Z_p[g]$-modules are finite, which gives the $p$-localized formula:
$$\frac{\order (\CH_{L,\varphi}^G)}{\order (\CH_{K,\varphi})} = 
\frac{\order \Coker(\J_{L/K})_\varphi}{\order \Ker(\J_{L/K})_\varphi} 
= \frac{\order \Coker(\jj_{L/K}^{})_\varphi }{(\CE_{K,\varphi} :
\CE_{K,\varphi} \cap \Norm_{L/K}(L^\times))} \times 
\frac{\order \Hom^2(G,\CE_{L,\varphi})}{\order \Hom^1(G,\CE_{L,\varphi})},$$

\noindent
where $\ds \frac{\order \Hom^2(G,\CE_{L,\varphi})}{\order \Hom^1(G,\CE_{L,\varphi})}$
is the Herbrand quotient of $\CE_{L,\varphi}$ we talked about, whose computation leads, 
for $\BE_L$, to the global value $\ffrac{1}{[L : K]}$ in the real case (cf. \cite{Lan1990}, 
\cite[Chap. IX, \S\S\,1,4]{Lan2000}). 

\begin{lemma}
The Herbrand quotient of $\CE_{L,\varphi}$ is trivial for all $\varphi \in \Phi_K \setminus \{1\}$.
\end{lemma}

\noindent{\bf Proof.}
We know that $(\BE_L \otimes \Q) \oplus \Q$ is the regular representation 
$\Q[G \times g]$ (see, e.g., \cite[Theorem I.3.7]{Gra2005});
so there exists a ``Minkowski unit'' $\varepsilon$ such that the $\Z[G \times g]$-module
generated by $\varepsilon$ is of finite index in $\BE_L$ that
one may choose prime to $p$; so $\CE_L$ is such that  $\CE_L \oplus \Z_p 
\simeq \Z_p[G][g]$ as $g$-modules. 
Thus, $(\CE_L \oplus \Z_p)_\varphi = \CE_{L,\varphi} \simeq
\Z_p[\mu_{g_\varphi}] [G]$ for $\varphi \ne 1$; whence the result.
\qed

\medskip
So, $\order (\CH_{L,\varphi}^G) =  \order (\CH_{K,\varphi}) \times
\ds \frac{\order \Coker(\jj_{L/K}^{})_\varphi }{(\CE_{K,\varphi} : \CE_{K,\varphi} 
\cap \Norm_{L/K}(L^\times))}$, for $\varphi \ne 1$, 
and $\CH_{L,1}^G = \CH_{K,1} = 1$, which completes the proof of the Theorem.
\qed

\medskip
\subsubsection{The main exact sequences associated to the filtration}\label{ssect}
We give, without proofs, the exact sequences leading to the formulas \eqref{filtration}
giving the orders of $\CH_L^{i+1}/\CH_L^i := (\CH_L/\CH_L^i)^G$ for $i \geq 1$,
and justifying the $p$-localizations of the formulas (see \cite[Section 3]{Gra2017} 
for the details). Moreover, we will use essentially the case of $(\CH_L/\CH_L^G)^G$.

\smallskip
Let ${\mathcal H}$ be a sub-$G$-module of $\CH_L$. Put 
$\wt {\mathcal H} = \{h \in \CH_L, \ \, h^{1-\sigma} \in {\mathcal H}\}$; so
$\big (\CH_L / {\mathcal H} \big)^G = \wt {\mathcal H} / {\mathcal H}$. 
We have the exact sequences:
$$1\too \CH_L^{G} \tooo  \wt {\mathcal H} \mathop{\tooo}^{1-\sigma} 
(\wt {\mathcal H})^{1-\sigma} \too 1  $$
$$1\too {}_{\Norm}{\mathcal H} \tooo  {\mathcal H}
\mathop{\tooo}^{\Norm_{L/K}} \Norm_{L/K} ({\mathcal H}) \too 1.$$

\noindent
Let ${\mathcal I} \subset I_L$ be such that 
${\mathcal I} P_L /P_L = {\mathcal H}$; so:
$$\Norm_{L/K}({\mathcal H}) = \Norm_{L/K}({\mathcal I}) P_K / P_K. $$
Let $\Lbda :=  \{x \in K^\times, \, (x) \in \Norm_{L/K}({\mathcal I})\}$;
the fundamental exact sequence is then:
$$1\to \big( \BE_K \Norm_{L/K}(L^\times) \big) \cap \Lbda
\too \Lbda \mathop{\too}^\varphi  {}_{\Norm}{\mathcal H} \big / 
(\wt {\mathcal H})^{1-\sigma} \to 1,$$

\noindent
where, for all $x \in \Lbda$, $\varphi(x) = {\mathfrak A} P_L \pmod{(\wt {\mathcal H})^{1-\sigma}}$, 
for any ${\mathfrak A} \in {\mathcal I}$ such that $\Norm_{L/K}({\mathfrak A}) = (x)$.
This exact sequence is $p$-localizable since $\Lbda$, containing $\BE_K$, is a
sub-$\Z$-module of finite type of $K^\times$. We deduce from the above:
\begin{equation}
(\wt \CH : \CH) = \ds
\frac{\order \CH_L^G \cdot \order (\wt \CH)^{1-\sigma}}
{\order \Norm_{L/K}(\CH) \cdot \order  {}_{\Norm}\CH} 
= \frac{\order \CH_L^G}
{\order \Norm_{L/K}(\CH) \cdot  ( {}_{\Norm}\CH: (\wt \CH)^{1-\sigma})} ;
\end{equation}

\noindent
thus $\order \big(\CH_L/ {\mathcal H} \big)^G 
=\ds \frac{\order  \CH_L^G}{\order \Norm_{L/K}({\mathcal H}) \cdot 
(\Lbda E_K \Norm_{L/K}(L^\times) : E_K \Norm_{L/K}(L^\times))}$.

\medskip
Then Chevalley--Herbrand formula gives the final $p$-localized result.

\begin{corollary} \label{filtra}
Let $K$ be a cyclic real field of prime-to-$p$ degree and let
$L \subset K(\mu_\ell^{})$, $\ell \equiv 1 \pmod {2 p^N}$, inert in $K$, with 
$[L : K] = p^n$, $n \in [1,N]$. Then for all $\varphi  \in \Phi_K$, the filtrations 
defined by $\CH^{i+1}_{L,\varphi}/\CH^i_{L,\varphi} := 
\big (\CH_{L,\varphi}/\CH^i_{L,\varphi}\big)^G$, fulfill
the higer rank Chevalley--Herbrand formulas:
\[\order \CH_{L,\varphi}^G = \order \CH_{K,\varphi} \ \ \hbox{and}\ \ 
\order \big(\CH^{i+1}_{L,\varphi}/\CH^i_{L,\varphi}\big) 
= \ds \frac{\order \CH_{K,\varphi}}{\order \Norm_{L/K}(\CH^i_{L,\varphi})}
\ \hbox{ for all $i\geq 1$}. \] 
\end{corollary}

\noindent{\bf Proof.} 
We then have a unique place $(\ell)$ of $K$, totally
ramified in $L/K$. So any $x \in \Lbda$ being norm of an ideal
(that we may choose prime to $(\ell)$), it is local norm at any place
distinct from $(\ell)$; then the product formula of the Hasse norm symbols
\cite[Theorem II.3.4.1]{Gra2005} gives that $x$ is everywhere local norm, hence 
global norm (Hasse's norm theorem). 

\smallskip
This applies to 
the group of units for which $\CE_{K,\varphi} \subset \Norm_{L/K}(L^\times)$,
whence $(\CE_{K,\varphi} : \CE_{K,\varphi} \cap \Norm_{L/K}(L^\times)) = 1$ (recall
that under the unicity of the ramified prime ideal in $L/K$,
$\Coker(\jj_{L/K}^{})_\varphi  = 1$ for $\varphi \ne 1$; see Remark \ref{coker0}).
\qed

\medskip
Of course, this does not mean $\CH_{L,\varphi}^G = \J_{L/K} (\CH_{K,\varphi})$
since there is most often capitulation of classes; this expresses
the subtlety of Chevalley--Herbrand formula for which we shall give another
description of $\CH_{L,\varphi}^G$ likely to involve the kernel of capitulation.

\section{Exact sequence of capitulation}\label{cap}
We consider a real cyclic extension $K/\Q$ of prime-to-$p$ degree and
$L=KL_0$ with $L_0/\Q$ cyclic of degree $p^n$, $n \geq 1$.
As for the case $L \subset K(\mu_\ell^{})$, we assume to simplify that any ramified prime
is totally ramified in $L/K$. Put $G := \Gal(L/K) =: \langle \sigma \rangle$.

\smallskip
Let ${\mathfrak A}P_{L} \in \BH_{L}$ be a class invariant under $G$; 
thus ${\mathfrak A}^{1-\sigma} = \alpha P_{L}$,
$\alpha \in L^\times$, and $\Norm_{L/K}(\alpha)$ is a unit $\varepsilon \in 
\BE_K \cap \Norm_{L/K}(L^\times)$; if $\alpha' = \alpha \eta$, $\eta \in \BE_{L}$,
is another generator of ${\mathfrak A}^{1-\sigma}$, then $\Norm_{L/K}(\alpha') = 
\varepsilon \Norm_{L/K}(\eta)$. This defines the map:
\[ \BH_{L}^G \too \BE_K \cap \Norm_{L/K}(L^\times)/\Norm_{L/K}(\BE_{L}), \]
which associates with ${\mathfrak A}P_{L} \in \BH_{L}^{G}$ the class of the unit
$\varepsilon = \Norm_{L/K}(\alpha)$.

\begin{theorem} \label{suitecap}
We have, for all $\varphi \in {\bf \Phi}_K$, the exact sequences:
$$1 \to \J_{L/K}(\CH_{K,\varphi}) \cdot \CH_{L,\varphi}^\ram  \to \CH_{L,\varphi}^{G} 
\to \CE_{K,\varphi} \cap \Norm_{L/K}(L^\times)/\Norm_{L/K}(\CE_{L,\varphi}) \to 1, $$
where $\CH_{L}^\ram \subseteq \CH_{L}^G$ is generated by the classes of the ramified 
prime ideals.
\end{theorem}

\noindent{\bf Proof.}
We shall establish the global exact sequence:
\[1 \to \J_{L/K}(\BH_K) \cdot \BH_{L}^\ram  \to \BH_{L}^{G} \to
\BE_K \cap \Norm_{L/K}(L^\times)/\Norm_{L/K}(\BE_{L}) \to 1.\]

(i) Image. Let $\varepsilon \in \BE_K \cap \Norm_{L/K}(L^\times)$; put
$\varepsilon = \Norm_{L/K}(\alpha)$, then the ideal $\alpha P_{L}$ being of norm $1$
is of the form ${\mathfrak A}^{1-\sigma}$, ${\mathfrak A} \in I_{L}$ (Lemma \ref{lemma0}\,(iii)), 
and its class is invariant giving the pre-image.

\smallskip
(ii) Kernel. Suppose that the image of the invariant class ${\mathfrak A} P_{L}$
is $\varepsilon = \Norm_{L/K}(\eta)$, $\eta \in \BE_{L}$; then 
$\Norm_{L/K}(\alpha) = \Norm_{L/K}(\eta)$ and $\Norm_{L/K}(\alpha \eta^{-1}) = 1$
giving $\alpha \eta^{-1} = \beta^{1-\sigma}$, $\beta \in L^\times$ (Hilbert's Theorem $90$) 
then $\alpha P_{L} = {\mathfrak A}^{1-\sigma} = (\beta P_{L})^{1-\sigma}$.
So, the class ${\mathfrak A} P_{L}$ is the class of the invariant ideal 
${\mathfrak A} (\beta)^{-1}$; but the group of invariants ideals is generated
by $\jj_{L/K}^{}(I_K)$ and the ramified primes of $L$ since ramification is total. 
Whence the result.
\qed

\medskip
\begin{corollary}\label{maincoro}
Let $K$ be a cyclic real field of prime-to-$p$ degree and let
$L \subset K(\mu_\ell^{})$, $\ell \equiv 1 \pmod {2 p^N}$, inert in $K$, with 
$[L : K] = p^n$, $n \in [1,N]$. Then:

\smallskip
(i) For all $\varphi  \in \Phi_K \setminus \{1\}$, we have 
$\order \CH_{L,\varphi}^G = \order \CH_{K,\varphi}$ (Corollary \ref{filtra})
and the exact sequences
$1 \to \J_{L/K}(\CH_{K,\varphi}) \to \CH_{L,\varphi}^G 
\to \CE_{K,\varphi} / \Norm_{L/K}(\CE_{L,\varphi}) \to 1$. 

\smallskip
(ii) The capitulation of $\CH_K$ in $L/K$ is equivalent to:
$$\CH_{L,\varphi}^G \simeq \CE_{K,\varphi} / \Norm_{L/K}(\CE_{L,\varphi}),
\ \hbox{for all $\varphi \in \Phi_K$}; $$
if so, $\CH_{L,\varphi}^G \simeq \Z_p[\mu_{g_\varphi}]/(p^{a_\varphi})\Z_p[\mu_{g_\varphi}]$,
$a_\varphi^{}$ being such that $p^{\rho_\varphi a_\varphi} = \order  \CH_{K,\varphi}$,
where $\rho_\varphi^{} = [\Q_p(\mu_{g_\varphi}) : \Q_p]$, and 
$\order \CH_{K,\varphi} = (\CE_{K,\varphi} : \Norm_{L/K}(\CE_{L,\varphi}))$.
\end{corollary}

\noindent{\bf Proof.}
Exact sequence in (i) comes from the fact that $\CH_{L}^\ram$, generated
by the prime ideal ${\mathfrak L} \mid {\mathfrak l}$ for the unique prime
ideal ${\mathfrak l} = (\ell)$ of $K$, is of character $\varphi=1$ and from
the fact that $\Coker(\jj_{L/K}^{})_\varphi = (\CE_{K,\varphi} :
\CE_{K,\varphi} \cap \Norm_{L/K}(L^\times)) = 1$ for $\varphi \ne 1$ (Remark \ref{coker0}).

\smallskip
Equivalence (ii) comes from the exact sequence. The last claims come from 
the monogenicity of $\CE_{K,\varphi}$ as $\Z_p[\mu_{g_\varphi}]$-module and
from the equality $\order \CH_{L,\varphi}^G = \order \CH_{K,\varphi}$.
\qed

\medskip
\subsection{Test of capitulation -- Numerical illustrations}
Recall that for $\ell \equiv 1 \pmod {p^N}$, inert in $K$, we consider
$L \subset K(\mu_\ell^{})$ of degree $p^n$ over $K$, $n \in [1,N]$.
In the computations $L$ is often denoted $K_n$ and $G$
is denoted $G_n$, etc.

\smallskip
To verify in practice some capitulations in $L/K$, we use the relation:
\[\J_{L/K}(\BH_K) = \J_{L/K}(\Norm_{L/K}(\BH_L)) = \Nu_{L/K}(\BH_L), \]
since $\Norm_{L/K} : \BH_L \to \BH_K$ is surjective ($L/K$ is totally ramified), 
and we compute the algebraic norm of the generators $h_i$ of $\BH_L$ given by PARI.
So, we obtain explicite relations $\Nu_{L/K}(h_i) = \prod_j h_j^{a_{i,j}}$ in 
$\BH_L$, the complete capitu\-lation being given by the identity $\Nu_{L/K}(\CH_L)
= \big\langle \, \ldots, \prod_j h_j^{a_{i,j}}, \ldots\,\big \rangle_{\!i} \otimes \Z_p = 1$
and incomplete capitulations are deduced from the matrices $(a_{i,j})$.

\smallskip
We will give numerical examples showing what happens, in the isomorphism (ii) of
the corollary since $\CE_{K,\varphi} / \Norm_{L/K}(\CE_{L,\varphi})$
is monogenic as Galois module, while $\CH_{K,\varphi}$ is not in general. 

\smallskip
(i) Let's give, first, a few words about the important PARI instruction ${\sf bnfisprincipal}$,
of constant use in the computations, to prove capitulation of a class ${\mathfrak a} P_K$
in $L$, whence that ${\mathfrak a} P_L$ is principal; it is described as follows by \cite{Pari2013}:

\ft\begin{verbatim}
bnfisprincipal(bnf,x,{flag=1}):bnf being output by bnfinit 
(with flag<=2), gives [v,alpha], where v is the vector of 
exponents on the class group generators and alpha is the 
generator of the resulting principal ideal. In particular 
x is principal if and only if v is the zero vector. 
\end{verbatim}\ns

Thus, the most important output is the vector of exponents like ${\sf [0, 0]}$ 
meaning total capitulation of the selected $p$-class. Nevertheless, these vectors 
may be only $0$ modulo the order of the $p$-classes.

\smallskip
(ii) Let's give an example for $K$ cubic and $p=2$.
\begin{example} \label{example0}
{\rm We consider a cyclic cubic field with $p=2$,  for which
$\CH_K \simeq \Z/4\Z \times \Z/4\Z$, $\CH_{K_1} \simeq \Z/24\Z \times 
\Z/8\Z \times \Z/2\Z \times \Z/2\Z$ and $\CH_{K_2} \simeq \Z/16\Z \times 
\Z/16\Z \times \Z/2\Z \times \Z/2\Z$, but there is complete capitulation of 
$\CH_K$ in $K_2$ for $\ell = 449$ inert in $K$ (${\sf CK, CK1, CK2}$
represent $\CH_K$, $\CH_{K_1}$, $\CH_{K_2}$, respectively):

\ft\begin{verbatim}
conductor f=2817 PK=x^3 - 939*x + 6886 CK=[12,4]
ell=449  N=2  Nn=2  n=1  CK=[12,4]  CK1=[24,8,2,2]
norm in K1/K of the component 1 of CK1: [12,0,0,1]~
norm in K1/K of the component 2 of CK1: [12,0,1,0]~
norm in K1/K of the component 3 of CK1: [0,0,0,0]~
norm in K1/K of the component 4 of CK1: [0,0,0,0]~

ell=449  N=2  Nn=2  n=2  CK=[12,4]  CK2=[48,16,2,2] 
norm in K2/K of the component 1 of CK2: [0,0,0,0]~
norm in K2/K of the component 2 of CK2: [0,0,0,0]~
norm in K2/K of the component 3 of CK2: [0,0,0,0]~
norm in K2/K of the component 4 of CK2: [0,0,0,0]~
\end{verbatim}\ns}
\end{example}

(iii) Let's give another example for $K$ quadratic and $p=3$.
\begin{example} \label{example1}
{\rm Let $m = 32009$, $K=\Q(\sqrt m)$ for which $\CH_K \simeq \Z/3\Z \times \Z/3\Z$. 
Take $\ell = 19$ (inert in $K$). The general Program \ref{programquad} gives an 
incomplete capitulation in $K_1$, then a total capitulation in $K_2$ (in this last data 
for $n=2$, we give the $18$ integer coefficients of a generator of the ideal on the 
integral basis computed by PARI):

\ft\begin{verbatim}
PK=x^2 - 32009 CK=[3,3]
ell=19  N=2  Nn=2  n=1  CK=[3,3]  CK1=[9,3]
norm in K1/K of the component 1 of CK1: [3,0]~
norm in K1/K of the component 2 of CK1: [0,0]~

ell=19  N=2  Nn=2  n=2  CK=[3,3]  CK2=[9,3]
norm in K2/K of the component 1 of CK2: [0,0]~
norm in K2/K of the component 2 of CK2: [0,0]~
[[0,0]~,
[1439216371631838382,473754414131112320,454228965737496519,
-536418025036156085,919689041243214339,207983767848706102,
-1036595574155274193,128338672307382267,840355575133838069,
14736364857686206,6548993298829283,8896045582984518,
3839859983910278,7389856141626720,10096758515087182,
-895338434531174,-1442732225425698,8004210362478777]~]
\end{verbatim}\ns

In $K_1/K$, the exact sequence looks like:
\[1 \to \J_{K_1/K}(\CH_K) \simeq \Z/3\Z \too \CH_{K_1}^{G_1} \too 
\CE_K / \Norm_{K_1/K}(\CE_{K_1}) \simeq \Z/3\Z \to 1,\]
the structure of $\CH_{K_1}^{G_1}$ being a priori unknown.
A direct computation shows that $\CH_{K_1}^{G_1} \simeq
\Z/3\Z \times \Z/3\Z$, but it is not $\J_{K_1/K}(\CH_K)$
since $\order \Ker(\J_{K_1/K}) = 3$.

\smallskip
In $K_2/K$, the exact sequence becomes the isomorphism:
\[\CH_{K_2}^{G_2} \simeq \CE_K/ \Norm_{K_2/K}(\CE_{K_2}) \simeq \Z/9\Z, \] 
since $\J_{K_2/K}(\CH_K)=1$ and $\CE_K \simeq \Z_3$.
So, we intend to find a generator of $\CH_{K_2}^{G_2}$.
Taking the class of order $9$ (first component of $\CH_{K_2}$ given by the instruction 
${\sf A0=Kn.clgp[3][1]}$), we compute its conjugate by the automorphism 
${\sf S}$ of order $9$:

${\sf B0=nfgaloisapply(Kn,S,A0)}$, then ${\sf C0=idealpow(Kn,B0,8)}$ and 

${\sf R=idealmul(Kn,A0,C0)}$ for which we apply the test:

${\sf U=bnfisprincipal(Kn,R)}$,

\noindent
giving a principal integer with huge integer coefficients:

\ft\begin{verbatim}
A0=Kn.clgp[3][1];B0=nfgaloisapply(Kn,S,A0);C0=idealpow(Kn,B0,8);
R=idealmul(Kn,A0,C0);U=bnfisprincipal(Kn,R);print(U)
[[0, 0],
[-15352694895259448716005913288179,-4937712840022370286191878614596, 
10234788031577460568990927971879,9473644150178364411380147768767, 
919093855688240643550377510520,4933150036914472598668945255159, 
-8186088867315265238068365774860,-12462519184163404099427848753248, 
-10116044365250124206306632744945,-32490191902858719490198341631, 
-25043275245516486863055415213,150430762033938424462889500018, 
83650512449675221885273689474,-243856275294066992198407658217, 
229543001775729020765314244564,-195555194837852410495787525040, 
140098686428490793581998673118,26378140206989683079187437611]~]
\end{verbatim}\ns

This confirms that $\CH_{K_2}^{G_2}$ is cyclic of order $9$ (cf.
Corollary \ref{maincoro}\,(ii)).

\smallskip
This phenomenon is general (when the capitulation is complete), for a $p$-class 
group of the form (say, for $K$ quadratic)
$\CH_K \simeq \Z/p^{a_1}\Z \times \cdots \times \Z/p^{a_r}\Z$ and 
gives, at a layer $n \geq a := a_1+ \cdots + a_r$, the isomorphism:
\[\CH_{K_n}^{G_n} \simeq \CE_K / \Norm_{K_n/K}(\CE_{K_n}) \simeq \Z/p^a\Z, \]
but $\order \CH_{K_n}^{G_n} = \order \CH_K$.
This is typical of the Main Conjecture philosophy and will be enforced by the analytic 
framework, recalled in the next subsection, showing that non-cyclic structures of the base 
field $K$ (i.e., that of $\CH_{K,\varphi}$) leads to canonical ones by extension in $L$
(i.e., as quotients of $\CE_{K,\varphi} \simeq \Z_p[\mu_{g_\varphi}^{}]$), 
whatever the $p$-rank and the exponent of $\CH_K$.

\smallskip
It seems that this ``monogenicity'', by suitable cyclic $p$-extensions, has not been
remarked in the literature. Unfortunately, proof of capitulations are perhaps 
out of reach despite their obviousness in the practice.}
\end{example}

\subsection{The stability as sufficient condition of capitulation}
Now, we give a sufficient condition of capitulation (see the comments given
in Remark \ref{stab}):

\begin{theorem}\label{criterion}
Consider a prime $\ell \equiv 1 \pmod {2 p^N}$, inert in the cyclic real field $K$, 
and $K_n \subset K(\mu_\ell^{})$, of degree $p^n$ over $K$, $n \in [1,N]$. 

\noindent
Then $\order \CH_{K_n} = \order \CH_K$ for all $n$ if and only if 
$\order \CH_{K_1} = \order \CH_K$. If this criterion applies, then 
$\CH_{K_n}^{G_n} = \CH_{K_n}$ for all $n$, 
$\Ker(\J_{K_n/K}) = \Norm_{K_n/K}(\CH_{K_n}[p^n])$, and if $p^e$ is the 
exponent of $\CH_K$, then $\CH_K$ capitulates in $K_e$
(assuming $N \geq e$).
\end{theorem}

\noindent{\bf Proof.} 
Consider $\Gal(K_n/K_1) = G_n^p$. Then we have the 
Chevalley--Herbrand formulas
$\order \CH_{K_n}^{G_n} = \order \CH_K$ and $\order \CH_{K_n}^{G_n^p} = \order \CH_{K_1}$. 
But $\CH_{K_n}^{G_n} \subseteq \CH_{K_n}^{G_n^p}$; then under the condition
$\order \CH_{K_1} = \order \CH_K$, we get $\CH_{K_n}^{G_n} = \CH_{K_n}^{G_n^p}$,
equivalent to 
$$\CH_{K_n}^{1-\sigma_n} = \CH_{K_n}^{1-\sigma_n^p} = 
\CH_{K_n}^{(1-\sigma_n) \,\cdot\, \theta}, $$
where $\theta = 1+\sigma_n+ \cdots + \sigma_n^{p-1} \in (p, 1-\sigma_n)$,
a maximal ideal of $\Z_p[G_n]$ since $\Z_p[G_n]/(p, 1-\sigma_n) \simeq \F_p$; 
so $\CH_{K_n}^{1-\sigma_n} = 1$, thus $\CH_{K_n}=\CH_{K_n}^{G_n}$ for all
$n \in [1, N]$. 

\smallskip
Reciprocal is trivial.

\smallskip
From $\Norm_{K_n/K}(\CH_{K_n}) = \CH_K$, $\CH_{K_n}=\CH_{K_n}^{G_n}$
and $\J_{K_n/K} \circ \Norm_{K_n/K} = \nu_{K_n/K}$, one gets
$\J_{K_n/K}(\CH_K) = \J_{K_n/K}(\Norm_{K_n/K}(\CH_{K_n})) = 
{\Nu_{K_n/K}}(\CH_{K_n} )= \CH_{K_n}^{p^n}$. 
Let $c \in \Ker(\J_{K_n/K})$ and put $c = \Norm_{K_n/K}(C)$, $C \in \CH_{K_n}$; so 
$1 = \J_{K_n/K}(c) = \J_{K_n/K}(\Norm_{K_n/K}(C))=C^{p^n}$, and 
$\Ker(\J_{K_n/K}) \subseteq \Norm_{K_n/K}(\CH_{K_n}[p^n])$.

\smallskip
Reciprocally, if $c=\Norm_{K_n/K}(C)$, $C^{p^n}=1$, then:
$$\J_{K_n/K}(c)=\J_n(\Norm_{K_n/K}(C)) = C^{p^n}=1;$$ 
whence:
$$\Ker(\J_{K_n/K}) = \Norm_{K_n/K}(\CH_{K_n}[p^n]) \subseteq \CH_K[p^n].$$ 

For $n = e$, one obtains the capitulation of $\CH_K$ in $K_e$.
\qed

\medskip
\begin{remarks}{\rm 
(i) Since all the relations and exact sequences defining the filtration $p$-localize,
the stability relation $\order \CH_{K_1,\varphi} = \order \CH_{K,\varphi}$ 
implies the capitulation of $\CH_{K,\varphi}$ in $K_e$. 

\smallskip
(ii) From formula of Corollary \ref{filtra}, the assumption $\order \CH_{K_1} = 
\order \CH_K$ is equivalent to $\ffrac{\order \CH_K}{\order \Norm_{K_1/K}(\CH_{K_1}^{G_1})} = 1$; 
indeed, if $\order \CH_{K_1} = \order \CH_K$, then from Theorem \ref{criterion}, 
$\CH_{K_1} = \CH_{K_1}^{G_1}$, so $\Norm_{K_1/K}(\CH_{K_1}^{G_1}) = 
\Norm_{K_1/K}(\CH_{K_1}) = \CH_K$.

\smallskip
If $\ffrac{\order \CH_K}{\order \Norm_{K_1/K}(\CH_{K_1}^{G_1})} = 1$, then
the filtration stops and $\CH_{K_1}^{G_1} = \CH_{K_1}$, whence 
$\order \CH_{K_1} = \order \CH_{K_1}^{G_1} = \order \CH_K$.

\smallskip
(iii) The same criterion holds if one replaces $K$ by $K_{n_0}$ for some $n_0 \geq 1$,
under the condition $N \geq n_0 + e$; the fact that $K_{n_0}/\Q$ is not of
prime-to-$p$ degree does not matter (proof of Theorem \ref{criterion} does
not need this assumption and requires that $L/K_{n_0}$ be totally ramified at a 
unique place; the notation $\CH_{K_n,\varphi}$ is still relative to characters 
of $K$, which makes sense since $\Gal(K_n/\Q) \simeq G_n \times g$).}
\end{remarks}

\section{Crucial link between Capitulation and Main Conjecture}

\subsection{Analytics -- The group of cyclotomic units}
This aspect being very classical, we just recall the main definitions and 
needed results. For the main definitions and properties of cyclotomic units, see 
\cite[\S\,8\ (1)]{Leo1954}, \cite{Leo1962} or \cite[Chap. 8]{Was1997}.

\begin{definition}
{\rm (i) Let $\chi \in \X$ even of conductor $f_\chi$; we define the ``cyclotomic numbers''
$\theta_\chi := \prod_{a \in A_\chi} (\zeta_{2f_\chi}^a - \zeta_{2f_\chi}^{-a})$, with
$\zeta_{2f_\chi} := \exp \big(\frac{i \pi}{f_\chi} \big)$, where $A_\chi$ is a half-system 
of representatives of $\Gal(\Q(\mu_{f_\chi}^{})/K_\chi )$ in $(\Z/ f_\chi \Z)^\times$.

\smallskip
(ii) Let $K$ be a real abelian field and let $\BF_K$ be the intersection with $\BE_K$
of the multiplicative group generated by the conjugates of $\theta_\chi$, for all 
$\chi \in \X_K$. }
\end{definition}

Recall that $\theta_\chi^2 \in K_\chi$ and that $\ffrac{\theta'_\chi}{\theta_\chi} \in \BE_{K_\chi}$
for any conjugate $\theta'_\chi$ of $\theta_\chi$. If $f_\chi$ is not a prime power, 
$\theta_\chi \in \BE_{K_\chi}$. Since we will consider $\varphi$-components of the
$\theta_\chi$, for $\varphi \ne 1$, one gets always units of $K_\chi$.

\smallskip
These units lead to an analytic computation of $\order \BH_\chi$, $\chi \in \X$ 
even, $\chi \ne 1$ (using Theorem \ref{chiformula}). 
One obtains, in the semi-simple case:

\begin{theorem}
For all $\chi \in \X$, even of prime-to-$p$ order, $\order \CH_\chi = (\CE_\chi : \CF_\chi)$.
\end{theorem}

The philosophy of the abelian Main Conjecture is to ask if the analogous relations
$\order \CH_\varphi = (\CE_\varphi : \CF_\varphi)$ exist or not, since we only know that:
\begin{equation}\label{prodchiformula}
\order \CH_\chi = \prd_{\varphi \mid \chi} \order \CH_\varphi =
(\CE_\chi : \CF_\chi) = \prd_{\varphi \mid \chi}(\CE_\varphi : \CF_\varphi). 
\end{equation}

\subsection{Norm properties of cyclotomic units}

We mention, first, the classical norm property of cyclotomic units that are given
in many books and articles, but are crucial for our purpose:

\begin{proposition}\label{cycloformula}
Let $f >1$ and then let $m \mid f$, with $m>1$, be any modulus; let 
$\Q^m := \Q(\zeta_m) \subseteq \Q^f  := \Q(\zeta_f)$ be the corresponding 
cyclotomic fields with $\zeta_t := \exp \big(\frac{2 i \pi}{t}  \big)$, for all $t \geq 1$.  
Put $\eta^{}_{\Q^f} := 1- \zeta_f$, $\eta^{}_{\Q^m} := 1- \zeta_m$; we have
$\Norm_{\Q^f/\Q^m} (\eta^{}_{\Q^f}) = 
 \eta_{\Q^m}^\Omega$, with $\Omega = {\prod_{\ell \mid f,\  \ell \nmid m}
\big(1-\big(\frac{\Q^m}{\ell} \big)^{-1}\big)}$,
where $\big(\frac{\Q^m}{\ell} \big) \in \Gal(\Q^m/\Q)$ denotes the Frobenius 
(or Artin) automorphism of the prime number $\ell \nmid m$, that is to say 
such that $\zeta_m \mapsto \zeta_m^\ell$.
\end{proposition}

\noindent{\bf Proof.}
To simplify, denote by $\tau^{}_a$, $a$ prime to $f$,
the Artin automorphism $\big(\frac{\Q^f}{a} \big)$ defined by
$\zeta_f \mapsto \zeta_f^a$, 
then put $\eta^{}_{\Q^f} =: \eta^{}_f$, $\eta^{}_{\Q^m} =: \eta^{}_m$.

\smallskip
We consider, by induction, the case $f=\ell \cdot m$,
with $\ell$ prime and examine the two cases $\ell \nmid m$ and $\ell \mid m$.
We have $\Norm_{\Q^f/\Q^m}(\eta^{}_f) = \prd_{a} \eta_f^{\tau_a}$
where $a$ runs trough the integers $a \in [1, f]$ prime to $f$ and such that
$a \equiv 1 \pmod m$.

\smallskip
(i) Case $\ell \nmid m$.
Put $a=1 + \lambda \cdot m$,  $\lambda \in [0, \ell-1]$, but
we must exclude a unique $\lambda^* \in [0, \ell-1]$ such that
$1 + \lambda^*\cdot m \equiv 0 \pmod \ell$; put $1+ \lambda^*m = \mu \ell$. 
Thus: 
$$\hspace{-3cm}\Norm_{\Q^f/\Q^m}(\eta^{}_f) = 
\prd_{\lambda \in [0, \ell-1],\, \lambda \ne \lambda^*} (1-\zeta_f^{1+\lambda m}) $$
$$\hspace{2.0cm} = \frac{\prod_{\lambda \in [0, \ell-1]}(1- \zeta_f \zeta_\ell^{\lambda})}
{1- \zeta_f ^{\mu \ell}} = \frac{1 - \zeta_f^\ell}{1 - \zeta_m^\mu} 
= \frac{1 - \zeta_m}{1 - \zeta_m^\mu}.$$

Since $\mu \equiv \ell^{-1} \pmod m$, we get
$\Norm_{\Q^f/\Q^m}(\eta^{}_f) = \eta_m^{1 - \tau_\ell^{-1}}$.

\smallskip
(ii) If $\ell \mid m$, any $\lambda \in [0, \ell-1]$
is suitable, giving $\Norm_{\Q^f/\Q^m} (\eta^{}_f)=\eta^{}_m$. 
\qed

\medskip
\begin{corollary} \label{invertible}
Let $L/K/\Q$ be real abelian extensions where 
$L$ is of conductor $f$ and $K$ of conductor $m$; set $\eta^{}_L := 
\Norm_{\Q^f/L}(\eta^{}_f)$ and $\eta^{}_K := \Norm_{\Q^m/K}(\eta^{}_m)$. Then:
\[\Norm_{L/K} (\eta^{}_L) = \eta_K^\Omega,\ \,\hbox{with \,$\Omega = 
\prod_{\ell \mid f,\  \ell \nmid m}
\big(1-\big(\frac{K}{\ell} \big)^{-1}\big)$} .\]

\noindent
If moreover $K/\Q$ is a cyclic extension of prime-to-$p$ degree, 
and if all the primes $\ell \mid f$, $\ell \nmid m$, are inert in $K$, 
then $\Omega e_\varphi = \prod_{\ell \mid f,\, \ell \nmid m} 
\big(1-\big(\frac{K}{\ell} \big)^{-1}\big) e_\varphi$ is an invertible element 
of the algebra $\Z_p[g]e_\varphi \simeq \Z_p[\mu_{g_\varphi}^{}]$, for all 
$\varphi \in {\bf \Phi}_K \setminus \{1\}$.
In particular, $\Norm_{L/K}{\CF_{L,\varphi}} = \CF_{K,\varphi}$ for all
$\varphi \in {\bf \Phi}_K \setminus \{1\}$.
\end{corollary}

\noindent{\bf Proof.}
Indeed, $\tau^{}_{\ell,K} := \big(\frac{K}{\ell} \big)$ is a generator of $g = \Gal(K/\Q)$
since $\ell$ is inert in $K/\Q$; so for $\psi \mid \varphi \ne 1$, 
$\psi \big(1-\tau^{-1}_{\ell,K}\big) = 1 - \psi\big(\tau^{-1}_{\ell,K}\big)$ is a unit
of $\Q(\mu_{g_\varphi}^{})$ if $\psi\big (\tau^{-1}_{\ell,K}\big)$ 
is not of prime power order, otherwise, if $g_\varphi$ is a power of a prime $q$, 
then since $q \ne p$, $\big(1 - \psi\big (\tau^{-1}_{\ell,K}\big)\big)$ is a prime ideal  
above $q$ in $\Q(\mu_{g_\varphi}^{})$ and $1 - \psi\big(\tau^{-1}_{\ell,K}\big)$ is 
a $p$-adic unit. Whence the norm relation between the $p$-localized groups of 
cyclotomic units for $\varphi \ne 1$.\qed

\begin{remark}
{\rm The link with the Leopoldt definition of cyclotomic units is easy since we get
$\zeta_{2f} - \zeta_{2f}^{-1} = - \zeta_{2f}^{-1}(1 - \zeta_{2f}^2) =
- \zeta_{2f}^{-1}(1 - \zeta_{f}) = - \zeta_{2f}^{-1} \eta^{}_f$, which has no
consequence for Proposition \ref{cycloformula} and its Corollary since norms are taken 
over real fields $L$, $K$. Numerically, $\eta_L$ and $\eta_K$
must be replaced by suitable square roots due to the use of the half-system $A_\chi$.}
\end{remark}

\subsection{Final statement}
So, we can state and prove the main result involving the transfer map $\J_{L/K}$
and its $p$-localized images $\J_{L/K}(\CH_{K,\varphi})$:

\begin{theorem}\label{relfond}
Let $K/\Q$ be a real cyclic extension of prime-to-$p$ degree.
Let $\ell \equiv 1 \pmod {2 p^N}$, $N \geq 1$, and assume $\ell$ totally inert in $K$.
Let $L \subset K(\mu_\ell^{})$ of degree $p^n$ over $K$, $n \in [1,N]$,
and put $G := \Gal(L/K) =: \langle \sigma \rangle$.

\smallskip
(i) We have the relations (product of two integers):
\[ (\CE_{K,\varphi} : \CF_{K,\varphi}) = \big (\Norm_{L/K} (\CE_{L,\varphi}) : \CF_{K,\varphi}\big)
\times \frac{\order \CH_{K,\varphi}}{\order \J_{L/K}(\CH_{K,\varphi})},\ \,
\hbox{for all $\varphi \in \Phi_K$}. \]

(ii) If $\CH_{K,\varphi}$ capitulates in $L$, then $(\CE_{K,\varphi} : \CF_{K,\varphi}) 
\geq \order \CH_{K,\varphi}$. 

\smallskip
(iii) The Main Conjecture $\order \CH_{K,\varphi} = (\CE_{K,\varphi} : \CF_{K,\varphi})$
for all $\varphi \in \Phi_K$, holds under the existence, for each 
$\varphi \in \Phi_K \setminus \{1\}$, of an inert prime $\ell \equiv 1 \pmod {2p^N}$, 
$N$ large enough, such that $\CH_{K,\varphi}$ capitulates in $K(\mu_\ell^{})$.
\end{theorem}

\noindent{\bf Proof.}
From Corollary \ref{maincoro}\,(i) to Theorem \ref{suitecap}, we have, 
for all $\varphi \ne 1$ the exact sequences
$1 \to \J_{L/K}(\CH_{K,\varphi}) \to \CH_{L,\varphi}^G 
\to \CE_{K,\varphi} / \Norm_{L/K}(\CE_{L,\varphi}) \to 1$, 
and $\order \CH_{L,\varphi}^G = \order \CH_{K,\varphi}$ from Corollary \ref{filtra};
whence the relations:
$$\order \CH_{L,\varphi}^G = \order \CH_{K,\varphi} =  
(\CE_{K,\varphi}  : \Norm_{L/K}(\CE_{L,\varphi})) \times \order \J_{L/K}(\CH_{K,\varphi}).$$ 

From Corollary \ref{invertible}, $ \CF_{K,\varphi} = \Norm_{L/K} (\CF_{L,\varphi})$, 
whence the inclusions:
\[\CF_{K,\varphi} \subseteq \Norm_{L/K} (\CE_{L,\varphi}) \subseteq \CE_{K,\varphi}, \]  
where $(\CE_{K,\varphi} : \Norm_{L/K} (\CE_{L,\varphi})) = 
\ds \frac{\order \CH_{K,\varphi}}{\order \J_{L/K}(\CH_{K,\varphi})}$, proving the 
claims (i) and (ii). For (iii), formula \eqref{prodchiformula},
$\prd_{\varphi \in \Phi_K}(\CE_{K,\varphi} : \CF_{K,\varphi}) = 
\prd_{\varphi  \in \Phi_K} \order \CH_{K,\varphi}$, implies equalities for
all $\varphi  \in \Phi_K$.
\qed

\medskip
Recall that a
sufficient condition of capitulation (whence implying the Main Conjecture) is the stability of the 
$p$-class groups in the cyclic $p$-tower $\bigcup_{n \in [1,N]} K_n$ of $K(\mu_\ell^{})/K$, 
that is to say, the existence of $n_0$, $0 \leq n_0 \leq N-e(\varphi)$, such that 
$\order \CH_{K_{n_0+1},\varphi} = \order \CH_{K_{n_0},\varphi}$, where 
$p^{e(\varphi)}$ is the exponent of $\CH_{K,\varphi}$.
Let's note that, as soon as there is an incomplete capitulation of $\CH_{K,\varphi}$ in some 
$K_n/K$, $n \in [1, N]$, the index $(\CE_{K,\varphi} : \CF_{K,\varphi}) $ is non trivial. 
In practice, one obtains often the whole capitulation of $\CH_K$ using a single prime 
$\ell$ among, probably, infinitely many.

\section{Numerical experiments over cyclic cubic fields} \label{comput}

We consider cyclic cubic fields $K$ with $p \ne 3$. Let $L \subset K(\mu_\ell^{})$, 
$\ell \equiv 1 \pmod {2 p^N}$ inert in $K$, be a cyclic $p$-extension of 
degree $p^N$ of $K$. 

\smallskip
The following program gives, at the beginning, the complete list of cyclic cubic fields of 
conductor ${\sf f \in [bf, Bf]}$ and selects those having a suitable (non-trivial) $p$-class 
group to study the capitulation in $L/K$. The test is about the order of $\CH_K$, so 
various structures may occur. The fields $K_n \subseteq L$, of degree $p^n$, are given 
by means of the polynomial ${\sf P}$ of degree ${\sf 3 \cdot p^n}$, ${\sf n \in [1, Nn]}$, 
where ${\sf Nn \leq N}$, not too large, defines the layers in which the computations 
are done (regrettably the execution time becomes rapidly out of reach).

\smallskip
As the instruction ${\sf bnfinit(P,1)}$ takes huge time if the 
degree of ${\sf P}$ increases, we are limited to $p=2$ and possibly 
$p=5$ and $7$ (minimal prime giving two $p$-adic characters). 
The purpose being to suggest the randomness of the (very frequent) phenomenon of 
capitulation, we hope that these cases constitute a good heuristic.
The case $p=2$ is not specific for capitulation aspects and allows
more complex structures for $\CH_K$, which is crucial to understand
the process in the tower.

\subsection{General program for cyclic cubic fields}

In the following general program, one must precise the following data:

\smallskip
(i) The numbers ${\sf N}$ (and ${\sf Nn \leq N}$, the number of layers to be
tested by the program) to define the primes ${\sf ell}$, limited by the bound 
${\sf Bell}$, congruent to $1$ modulo ${\sf 2 p^{N}}$; it is possible, to prove 
capitulations at a larger layer, to take ${\sf N}$ large, but ${\sf Nn}$ very
small.

\smallskip
(ii) The bounds ${\sf bf, Bf}$ defining an interval for the conductors ${\sf f}$.

\smallskip
(iii) The positive numbers ${\sf vHK}$ and ${\sf vHKn}$, for the
instructions:

\smallskip
\centerline{${\sf valuation(HK,p)<vHK}$, ${\sf valuation(HKn,p)<vHKn}$,}

\smallskip\noindent
to get only interesting $p$-class groups for ${\sf K}$ and ${\sf Kn}$; note that from the 
Chevalley--Herbrand formula, $\order \CH_{K_n}$ (in ${\sf HKn}$) is a multiple 
of $\order \CH_K$  (in ${\sf HK}$), 
and one must take ${\sf vHK \leq vHKn}$.

\ft\begin{verbatim}
STUDY OF THE CAPITULATION OF HK IN Kn/K FOR CYCLIC CUBIC FIELDS 
{p=2;N=2;Nn=2;bf=7;Bf=5*10^3;vHK=4;vHKn=6;Bell=500;
\\List of cubic fields of any conductor f:
for(f=bf,Bf,h=valuation(f,3);if(h!=0 & h!=2,next);
F=f/3^h;if(core(F)!=F,next);
F=factor(F);Div=component(F,1);d=matsize(F)[1];for(j=1,d,D=Div[j];
if(Mod(D,3)!=1,break));for(b=1,sqrt(4*f/27),
if(h==2 & Mod(b,3)==0,next);A=4*f-27*b^2;
if(issquare(A,&a)==1,\\a and b such that f=(a^2+27b^2)/4
if(h==0,if(Mod(a,3)==1,a=-a);PK=x^3+x^2+(1-f)/3*x+(f*(a-3)+1)/27);
if(h==2,if(Mod(a,9)==3,a=-a);PK=x^3-f/3*x-f*a/27);
\\End of computation of PK defining K of conductor f.
K=bnfinit(PK,1);HK=K.no;\\Whole class number of K
\\Test on the order of the p-class group of K:
if(valuation(HK,p)<vHK,next);CK=K.clgp;\\Class group of K
\\Definition of the primes ell inert in K:
forprime(ell=1,Bell,if(Mod(ell-1,2*p^N)!=0 || Mod(f,ell)==0,next);
F=factor(PK+O(ell));if(matsize(F)[1]!=1,next);
\\Definitions of the fields Kn<L, computation of their class group:
for(n=1,Nn,QKn=polsubcyclo(ell,p^n);P=polcompositum(PK,QKn)[1];
Kn=bnfinit(P,1);
HKn=Kn.no;htame=HKn/p^valuation(HKn,p);\\Tame part of HKn
\\Test on the order of the p-class group of Kn:
if(valuation(HKn,p)<vHKn,break);
print();print("conductor f=",f," PK=",PK," CK=",CK[2]);
print("ell=",ell," N=",N," Nn=",Nn," n=",n);
CKn=Kn.clgp;print("CKn=",CKn[2]);
rKn=matsize(CKn[2])[2];\\rank of CKn
\\Calcul de Gal(Kn/K) and search of a generator S of order p^n:
G=nfgaloisconj(Kn);Id=G[1];for(k=2,3*p^n,Z=G[k];ks=1;while(Z!=Id,
Z=nfgaloisapply(Kn,G[k],Z);ks=ks+1);if(ks==p^n,S=G[k];break));
\\Computation of the algebraic norms of the generators of CKn:
for(j=1,rKn,A0=CKn[3][j];A=1;
for(t=1,p^n,As=nfgaloisapply(Kn,S,A);A=idealmul(Kn,A0,As));
A=idealpow(Kn,A,htame);\\Replace A by its p-component
\\Test of capitulation (incomplete or total):
X=bnfisprincipal(Kn,A)[1];print("norm in K",n,"/K 
of the component ",j," of CK",n,": ",X)))))))}
\end{verbatim}\ns

\subsection{Case of cyclic cubic fields and \texorpdfstring{$p=2$}{Lg}}
We give an excerpt of the various forms of examples, with the structure
$\CH_K \simeq \Z/2\Z \times \Z/2\Z$ in $K$ and $\CH_{K_1}$ of order 
at least $2^6$ (we indicate the nature of the capitulation at the end of the data):

\ft\begin{verbatim}
conductor f=1777 PK=x^3 + x^2 - 592*x + 724 CK=[4,4] 
ell=41  N=3  Nn=2  n=1  CK=[4,4]  CK1=[4,4,2,2] 
norm in K1/K of the component 1 of CK1: [0,0,0,1]~
norm in K1/K of the component 2 of CK1: [2,2,1,1]~
norm in K1/K of the component 3 of CK1: [0,0,0,0]~
norm in K1/K of the component 4 of CK1: [0,0,0,0]~
Incomplete capitulation in K1

ell=41  N=3  Nn=2  n=2  CK=[4,4]  CK2=[8,8,2,2]  
norm in K2/K of the component 1 of CK2: [0,0,0,0]~
norm in K2/K of the component 2 of CK2: [0,0,0,0]~
norm in K2/K of the component 3 of CK2: [0,0,0,0]~
norm in K2/K of the component 4 of CK2: [0,0,0,0]~
Complete capitulation in K2 without stabilization in K2/K1
_____
ell=113  N=3  Nn=2  n=1  CK=[4,4]  CK1=[4,4,2,2] 
norm in K1/K of the component 1 of CK1: [0,0,1,0]~
norm in K1/K of the component 2 of CK1: [0,0,0,1]~
norm in K1/K of the component 3 of CK1: [0,0,0,0]~
norm in K1/K of the component 4 of CK1: [0,0,0,0]~
Incomplete capitulation in K1

ell=113  N=3  Nn=2  n=2  CK=[4,4]  CK2=[8,8,2,2] 
norm in K2/K of the component 1 of CK2: [0,4,0,0]~
norm in K2/K of the component 2 of CK2: [4,4,0,0]~
norm in K2/K of the component 3 of CK2: [0,0,0,0]~
norm in K2/K of the component 4 of CK2: [0,0,0,0]~
Incomplete capitulation in K2
_____
ell=257  N=3  Nn=2  n=1 CK=[4,4]  CK1=[72,24] 
norm in K1/K of the component 1 of CK1: [54,12]~
norm in K1/K of the component 2 of CK1: [36,18]~
No capitulation in K1

ell=257  N=3  Nn=2  n=2  CK=[4,4]  CK2=[72,24] 
norm in K2/K of the component 1 of CK2: [36,0]~
norm in K2/K of the component 2 of CK2: [0,12]~
Incomplete capitulation in K2
Stability from K1--->capitulation in K3
_____
ell=337  N=3  Nn=2  n=1  CK=[4,4]  CK1=[4,4] 
norm in K1/K of the component 1 of CK1: [2,0]~
norm in K1/K of the component 2 of CK1: [0,2]~
Incomplete capitulation in K1
Stability from K--->capitulation in K2

ell=337 N=3 Nn=3 n=2  CK=[4,4]  CK2=[4,4] 
norm in K2/K of the component 1 of CK2: [0, 0]~
norm in K2/K of the component 2 of CK2: [0, 0]~
Complete capitulation in K2 as expected
_____
ell=2129  N=3  Nn=2  n=1  CK=[4,4]  CK1=[16,16,2,2] 
norm in K1/K of the component 1 of CK1: [8,4,0,1]~
norm in K1/K of the component 2 of CK1: [12,12,1,0]~
norm in K1/K of the component 3 of CK1: [8,8,0,0]~
norm in K1/K of the component 4 of CK1: [0,8,0,0]~
No capitulation in K1

ell=2129  N=3  Nn=2  n=2  CK=[4,4]  CK2=[16,16,2,2] 
norm in K2/K of the component 1 of CK2: [0,8,0,0]~
norm in K2/K of the component 2 of CK2: [8,8,0,0]~
norm in K2/K of the component 3 of CK2: [0,0,0,0]~
norm in K2/K of the component 4 of CK2: [0,0,0,0]~
Incomplete capitulation in K2
Stability from K1--->capitulation in K_3
-----------------------------------------------------
conductor f=2817 PK=x^3 - 939*x + 6886 CK=[12,4]
ell=449  N=2  Nn=2  n=1  CK=[4,4]  CK1=[24,8,2,2] 
norm in K1/K of the component 1 of CK1: [12,0,0,1]~
norm in K1/K of the component 2 of CK1: [12,0,1,0]~
norm in K1/K of the component 3 of CK1: [0,0,0,0]~
norm in K1/K of the component 4 of CK1: [0,0,0,0]~
Incomplete capitulation in K1

ell=449  N=2  Nn=2  n=2  CK=[4,4]  CK2=[48,16,2,2] 
norm in K2/K of the component 1 of CK2: [0,0,0,0]~
norm in K2/K of the component 2 of CK2: [0,0,0,0]~
norm in K2/K of the component 3 of CK2: [0,0,0,0]~
norm in K2/K of the component 4 of CK2: [0,0,0,0]~
Complete capitulation in K2
-----------------------------------------------------
conductor f=4297 PK=x^3 + x^2 - 1432*x + 20371 CK=[4,4]
ell=449  N=2  Nn=2  n=1  CK=[4,4]  CK1=[4,4,2,2] 
norm in K1/K of the component 1 of CK1: [0,2,1,1]~
norm in K1/K of the component 2 of CK1: [2,2,0,1]~
norm in K1/K of the component 3 of CK1: [0,0,0,0]~
norm in K1/K of the component 4 of CK1: [0,0,0,0]~
Incomplete capitulation in K1

ell=449  N=2  Nn=2  n=2  CK=[4,4]  CK2=[292,4,4,4] 
norm in K2/K of the component 1 of CK2: [146,0,2,0]~
norm in K2/K of the component 2 of CK2: [0,0,2,2]~
norm in K2/K of the component 3 of CK2: [146,0,2,0]~
norm in K2/K of the component 4 of CK2: [146,0,2,0]~
Incomplete capitulation in K2
-----------------------------------------------------
conductor f=5409 PK=x^3 - 1803*x + 29449 CK=[12,4]
ell=113  N=2  Nn=2  n=1  CK=[4,4]  CK1=[12,4,2,2] 
norm in K1/K of the component 1 of CK1: [6,2,0,1]~
norm in K1/K of the component 2 of CK1: [0,0,1,0]~
norm in K1/K of the component 3 of CK1: [0,0,0,0]~
norm in K1/K of the component 4 of CK1: [0,0,0,0]~
Incomplete capitulation in K1

ell=113  N=2  Nn=2  n=2  CK=[4,4]  CK2=[24,8,2,2,2,2] 
norm in K2/K of the component 1 of CK2: [0,4,0,0,0,0]~
norm in K2/K of the component 2 of CK2: [12,4,0,0,0,0]~
norm in K2/K of the component 3 of CK2: [0,0,0,0,0,0]~
norm in K2/K of the component 4 of CK2: [0,0,0,0,0,0]~
norm in K2/K of the component 5 of CK2: [0,0,0,0,0,0]~
norm in K2/K of the component 6 of CK2: [0,0,0,0,0,0]~
Incomplete capitulation in K2
\end{verbatim}\ns

\begin{remark}
{\rm For the first example above, the capitulation in $K_2$ is complete, even if 
the stability does not occur from the first layer; the step $n=1$ shows an incomplete 
capitulation giving, $\J_{K1/K}(\CH_K) \simeq \Z/2\Z \times \Z/2\Z$ 
(indeed the exponent of $\CH_K$ is $4$).

\smallskip
To be more convincing, let's give the coefficients, on the PARI integral basis, 
of the generators of the representative ideals of the classes $\Nu_{K_2/K}(h_i)$
for the four classes $h_i$ of orders $8$, $8$, $2$, $2$, respectively (and given
as usual by ${\sf CKn=Kn.clgp}$); for this, replace,  in the program, 
${\sf X=bnfisprincipal(Kn,A)[1]}$ by ${\sf X=bnfisprincipal(Kn,A)}$. 
One gets $12$ huge integer coefficients:

\ft\begin{verbatim}
conductor f=1777   PK=x^3 + x^2 - 592*x + 724   CK=[4,4] 
ell=41  N=3  Nn=2  n=2  CK=[4,4]  CK1=[4,4,2,2]  CK2=[8,8,2,2]

norm in K2/K of the component 1 of CK2:
[[0,0,0,0]~,
[-31780222254443,-12898232803596,15554698429537,-11030242667244,
699234644603,-1433180593820,-97846196830,480428807611,
70679128541,581754438178,701511836521,497443446811]~]

norm in K2/K of the component 2 of CK2:
[[0,0,0,0]~,
[-409139735188114,-166218206982845,200303000159397,-142025979393819,
10482585098180,32200771380552,-1471046798500,7172751349077,
1043934050162,-13082943399099,-15769847218663,-11179282571565]~]

norm in K2/K of the component 3 of CK2:
[[0,0,0,0]~,
[4595853941743,7574362186256,-7431095890343,3180376719682,
-878235409486,520990038484,351933447679,127583225152,
327914236819,-381696290156,-181901226812,-173412643330]~]

norm in K2/K of the component 4 of CK2:
[[0,0,0,0]~,
[206178918528161385818507,45009133745540603328818,
109639228343931043367320,18671957635985615261071,
13770653000372358842954,7958894412958725580875,
13213450239959129254028,3845666006771496793309,
2250298427236450403785,-1737437297938711409589,
4232280623726481124024,-720772717060054219812]~]
\end{verbatim}\ns}
\end{remark}

\subsection{Case of cyclic cubic fields and \texorpdfstring{$p=7$}{Lg}}
For $p=7$, due to the execution time, let's give some examples of the case
$n=1$ with $\CH_K$ of order $7$, then one case of order $7^2$, and $B_\ell=100$.

\smallskip
We obtain complete capitulations in $K_1$, except few cases;
we give an excerpt of some possibilities:

\ft\begin{verbatim}
conductor f=313 PK=x^3 + x^2 - 104*x + 371 CK=[7]
ell=29  N=1  Nn=1  n=1  CK=[7]  CK1=[7]
norm in K1/K of the component 1 of CK1: [0]~
Complete capitulation in K1
Stability from K
-----------------------------------------------------
conductor f=1261 PK=x^3 + x^2 - 420*x - 1728 CK=[21]
ell=43  N=1  Nn=1  n=1  CK=[7]  CK1=[21]
norm in K1/K of the component 1 of CK1: [0]~
Complete capitulation in K1
Stability from K
-----------------------------------------------------
conductor f=1567 PK=x^3 + x^2 - 522*x - 4759 CK=[7]
ell=29  N=1  Nn=1  n=1  CK=[7] CK1=[49]
norm in K1/K of the component 1 of CK1: [7]~
No capitulation in K1
-----------------------------------------------------
conductor f=8563  PK=x^3 + x^2 - 2854*x + 57721  CK=[49]
ell=71  N=1  Nn=1  n=1 CK=[49] CK1=[49] 
norm in K1/K of the component 1 of CK1: List([7])
Incomplete capitulation in K1
Stability from K but N too small (e=2)
\end{verbatim}\ns

\smallskip
The last case shows an incomplete capitulation giving $\J_{K1/K}(\CH_K) \simeq \Z/7\Z$.
Since $N=1$, there is no possible complete capitulation despite the stability from $K$.
The case of primes $\ell$ with $N=n=2$ seems out of reach.

\smallskip
Let's give the generator of the principal ideal obtained after capitulation (first example above):

\ft\begin{verbatim}
conductor f=313  PK=x^3 + x^2 - 104*x + 371 CK=[7]
ell=29  N=1  Nn=1  n=1  CK=[7]  CK1=[7]
norm in K1/K of the component 1 of CK1:
[[0]~,
[4529357,2479589,125622,-2879283,2922668,4270474,-6202812,
-107453,1865872,37436,-613198,1546287,-1637834,1355628,
1276626,886508,944469,-900999,474890,508450,962907]~]
\end{verbatim}\ns

\subsection{Application of the sufficient condition of capitulation}
We consider cyclic cubic fields and $p=2$. We only search
examples of primes $\ell$ giving the stability of the $2$-class groups in $K_1/K$,
so that the capitulation automatically applies in $L=K_N$ if $N$ is large enough; 
so we compute the $2$-valuation of $\order \CH_K$ (in ${\sf N=valuation(HK,p)}$), 
which is the minimal possible bound (${\sf N}$ must be even).
We illustrate the case $\CH_K \simeq \Z/2\Z \times \Z/2\Z$; a great lot of examples
are found giving capitulation in $K_1$, even with small $\ell$'s (one writes only 
${\sf f, PK, ell, CK, CK1}$; we put ${\sf ell = prime(t)}$, ${\sf t \in [2,nell]}$):

\ft\begin{verbatim}
{p=2;N=6;bf=7;Bf=10^6;nell=100;for(f=bf,Bf,h=valuation(f,3);
if(h!=0 & h!=2,next);F=f/3^h;if(core(F)!=F,next);F=factor(F);
Div=component(F,1);d=matsize(F)[1];
for(j=1,d,D=Div[j];if(Mod(D,3)!=1,break));
for(b=1,sqrt(4*f/27),if(h==2 & Mod(b,3)==0,next);
A=4*f-27*b^2;if(issquare(A,&a)==1,
if(h==0,if(Mod(a,3)==1,a=-a);PK=x^3+x^2+(1-f)/3*x+(f*(a-3)+1)/27);
if(h==2,if(Mod(a,9)==3,a=-a);PK=x^3-f/3*x-f*a/27);
K=bnfinit(PK,1);HK=K.no;v=valuation(HK,p);if(v!=N,next);
for(t=2,nell,ell=prime(t);if(Mod(ell-1,2*2^N)!=0,next);
if(Mod(f,ell)==0,next);F=factor(PK+O(ell));
if(matsize(F)[1]!=1,next);QK1=polsubcyclo(ell,p);
P=polcompositum(PK,QK1)[1];K1=bnfinit(P,1);
if(valuation(K1.no,p)==N,print("f=",f," PK=",PK," ell=",ell,
"  CK=",K.clgp[2]," CK1=",K1.clgp[2]);break)))))}
N=2
f=163  PK=x^3+x^2-54*x-169         ell=29  CK=[2,2]  CK1=[2,2]
f=277  PK=x^3+x^2-92*x+236         ell=5   CK=[2,2]  CK1=[2,2]
f=349  PK=x^3+x^2-116*x-517        ell=5   CK=[2,2]  CK1=[2,2]
f=397  PK=x^3+x^2-132*x-544        ell=5   CK=[2,2]  CK1=[2,2]
(...)
f=9709 PK=x^3+x^2-3236*x+21216     ell=5   CK=[6,6]  CK1=[6,6]
f=9721 PK=x^3+x^2-3240*x-39244     ell=5   CK=[2,2]  CK1=[2,2]
f=9891 PK=x^3-3297*x+70336         ell=29  CK=[6,6]  CK1=[6,6]
f=9961 PK=x^3+x^2-3320*x-74523     ell=5   CK=[6,2]  CK1=[6,2]
(...)
-----------------------------------------------------------------
N=4
f=1777 PK=x^3+x^2-592*x+724        ell=353 CK=[4,4]  CK1=[52,4]
f=2817 PK=x^3-939*x+6886           ell=97  CK=[12,4] CK1=[12,12]
f=4297 PK=x^3+x^2-1432*x+20371     ell=97  CK=[4,4]  CK1=[4,4]
f=5409 PK=x^3-1803*x+29449         ell=193 CK=[12,4] CK1=[12,4]
(...)
f=98479 PK=x^3+x^2-32826*x-1940401 ell=353 CK=[4,4]  CK1=[4,4]
f=98581 PK=x^3+x^2-32860*x+1453157 ell=193 CK=[12,4] CK1=[228,12]
f=99133 PK=x^3+x^2-33044*x+117491  ell=193 CK=[4,4]  CK1=[4,4]
f=100807 PK=x^3+x^2-33602*x+321089 ell=97  CK=[12,4] CK1=[12,4]
(...)
-----------------------------------------------------------------
N=6
f=10513 PK=x^3+x^2-3504*x-80989    ell=257 CK=[8,8]  CK1=[24,8]
f=48769 PK=x^3+x^2-16256*x-7225    ell=257 CK=[24,8] CK1=[24,24]
f=70897 PK=x^3+x^2-23632*x-1389056 ell=257 CK=[24,8] CK1=[72,24,3]
f=80947 PK=x^3+x^2-26982*x+1696889 ell=257 CK=[24,8] CK1=[168,24,3,3]
(...)
f=351063 PK=x^3-117021*x-15407765    ell=257 CK=[24,24,3]CK1=[72,24,3,3]
f=357229 PK=x^3+x^2-119076*x+15228540ell=257 CK=[8,8]    CK1=[24,8]
f=492517 PK=x^3+x^2-164172*x-24479919ell=257 CK=[24,8]   CK1=[168,168,3]
f=552763 PK=x^3+x^2-184254*x-27842877ell=257 CK=[24,8]   CK1=[24,24,3]
(...)
\end{verbatim}\ns

Some cases of more complex structures of $\CH_K$ do not give stabilization
at the first step; this clearly depends on the exponent of the 
class groups as we have explained; but this sufficient condition is not necessary 
and capitulation does appear at larger layers; we illustrate (using the general program) 
the cases $\CH_K \simeq (\Z/2\Z)^4$ or $\CH_K \simeq (\Z/4\Z)^2 \times (\Z/2\Z)^2$:

\ft\begin{verbatim}
conductor f=7687 PK=x^3 + x^2 - 2562*x - 48969 CK=[2,2,2,2]
ell=17  N=3  Nn=3  n=1  CK=[2,2,2,2]  CK1=[4,4,2,2]
norm in K1/K of the component 1 of CK1: [2,2,0,0]~
norm in K1/K of the component 2 of CK1: [0,2,0,0]~
norm in K1/K of the component 3 of CK1: [2,0,0,0]~
norm in K1/K of the component 4 of CK1: [0,2,0,0]~
Incomplete capitulation in K1

ell=17  N=3  Nn=3  n=2  CK=[2,2,2,2]  CK2=[4,4,2,2]
norm in K2/K of the component 1 of CK2: [0,0,0,0]~
norm in K2/K of the component 2 of CK2: [0,0,0,0]~
norm in K2/K of the component 3 of CK2: [0,0,0,0]~
norm in K2/K of the component 4 of CK2: [0,0,0,0]~
Complete capitulation in K2 with stability from K1
-----------------------------------------------------
conductor f=20887 PK=x^3 + x^2 - 6962*x - 225889 CK=[4,4,2,2]
ell=193  N=3  Nn=3  n=1  CK=[4,4,2,2]  CK1=[8,8,2,2]
norm in K1/K of the component 1 of CK1: [2,4,0,0]~
norm in K1/K of the component 2 of CK1: [4,6,0,0]~
norm in K1/K of the component 3 of CK1: [4,0,0,0]~
norm in K1/K of the component 4 of CK1: [4,4,0,0]~
Incomplete capitulation in K1

ell=193  N=3  Nn=3  n=2  CK=[4,4,2,2]  CK2=[8,8,2,2]
norm in K2/K of the component 1 of CK2: [4,0,0,0]~
norm in K2/K of the component 2 of CK2: [0,4,0,0]~
norm in K2/K of the component 3 of CK2: [0,0,0,0]~
norm in K2/K of the component 4 of CK2: [0,0,0,0]~
Incomplete capitulation in K2
Stability from K1--->Complete capitulation in K3
-----------------------------------------------------
conductor f=31923 PK=x^3 - 10641*x + 227008 CK=[6,2,2,2]
ell=97  N=3  Nn=3  n=1  CK=[6,2,2,2]  CK1=[12,4,2,2,2,2]
norm in K1/K of the component 1 of CK1: [0,0,1,0,0,0]~
norm in K1/K of the component 2 of CK1: [0,0,0,1,0,0]~
norm in K1/K of the component 3 of CK1: [0,0,0,0,0,0]~
norm in K1/K of the component 4 of CK1: [0,0,0,0,0,0]~
norm in K1/K of the component 5 of CK1: [6,0,0,0,0,0]~
norm in K1/K of the component 6 of CK1: [6,2,0,0,0,0]~
No capitulation in K1

ell=97  N=3  Nn=3  n=2  CK=[6,2,2,2]  CK2=[312,8,4,4,4,4]
norm in K2/K of the component 1 of CK2: [156,0,2,0,0,2]~
norm in K2/K of the component 2 of CK2: [156,4,0,0,2,2]~
norm in K2/K of the component 3 of CK2: [156,0,0,0,0,0]~
norm in K2/K of the component 4 of CK2: [156,4,0,0,0,0]~
norm in K2/K of the component 5 of CK2: [156,0,0,0,0,0]~
norm in K2/K of the component 6 of CK2: [156,0,0,0,0,0]~
No capitulation in K2
\end{verbatim}\ns

So we must try another $\ell$ for $f=31923$:

\ft\begin{verbatim}
conductor f=31923 PK=x^3 - 10641*x + 227008 CK=[6,2,2,2]
ell=257  N=3  Nn=2  n=1 CK=[6,2,2,2]  CK1=[18,6,2,2,2,2]
norm in K1/K of the component 1 of CK1: [0,0,0,0,0,0]~
norm in K1/K of the component 2 of CK1: [0,0,0,0,0,0]~
norm in K1/K of the component 3 of CK1: [0,3,0,0,0,1]~
norm in K1/K of the component 4 of CK1: [9,0,0,0,1,0]~
norm in K1/K of the component 5 of CK1: [0,0,0,0,0,0]~
norm in K1/K of the component 6 of CK1: [0,0,0,0,0,0]~
Incomplete capitulation in K1

ell=257  N=3  Nn=2  n=2  CK=[6,2,2,2]  CK2=[36,12,2,2,2,2]
norm in K2/K of the component 1 of CK2: [0,0,0,0,0,0]~
norm in K2/K of the component 2 of CK2: [0,0,0,0,0,0]~
norm in K2/K of the component 3 of CK2: [0,0,0,0,0,0]~
norm in K2/K of the component 4 of CK2: [0,0,0,0,0,0]~
norm in K2/K of the component 5 of CK2: [0,0,0,0,0,0]~
norm in K2/K of the component 6 of CK2: [0,0,0,0,0,0]~
Complete capitulation in K2
\end{verbatim}\ns

For this new $\ell$, the complete capitulation is obtained in $K_2$. We note that 
the $2$-rank of $\CH_K$ is $4$, which implies, from Corollary \ref{maincoro}\,(ii)
and $j = \exp(\frac{2 i \pi}{3})$, the isomorphism $\CH_{K_2}^{G_2} \simeq \Z_2[j]/(4)$ 
but never $\Z_2[j]/(2) \times \Z_2[j]/(2)$.

\subsection{Statistics}
Let's consider an example of cyclic cubic field $K$ such that, for $p=2$, $\CH_K \simeq 
\Z/4\Z \times \Z/4\Z \times\Z/2\Z \times\Z/2\Z$ (conductor $20887$); we give some heuristics 
about the capitulation of $\CH_K$ in $K_1$ and $K_2$ (necessarily incomplete in $K_1$), 
for all prime numbers $\ell \equiv 1 \pmod 8$.

\ft\begin{verbatim}
{p=2;Nn=2;PK=x^3+x^2-6962*x-225889;Bell=10^3;
K=bnfinit(PK,1);CK=K.clgp[2];print("PK=",PK," CK=",CK);
forprime(ell=5,Bell,N=valuation((ell-1)/2,p);
if(N<Nn+1,next);F=factor(PK+O(ell));if(matsize(F)[1]!=1,next);
for(n=1,Nn,print();print("ell=",ell," N=",N," Nn=",Nn," n=",n);
QKn=polsubcyclo(ell,p^n);P=polcompositum(PK,QKn)[1];
Kn=bnfinit(P,1);CKn=Kn.clgp;print("CK=",CK,"  CK",n,"=",CKn[2]);
htame=Kn.no/p^valuation(Kn.no,p);rKn=matsize(CKn[2])[2];
G=nfgaloisconj(Kn);Id=G[1];for(k=2,3*p^n,Z=G[k];ks=1;
while(Z!=Id,Z=nfgaloisapply(Kn,G[k],Z);ks=ks+1);
if(ks==p,S=G[k];break));for(j=1,rKn,A0=CKn[3][j];A=1;
for(t=1,p^n,As=nfgaloisapply(Kn,S,A);A=idealmul(Kn,A0,As));
A=idealpow(Kn,A,htame);X=bnfisprincipal(Kn,A)[1];
print("norm in K",n,"/K of the component ",j," of CK",n,": ",X))))}
\end{verbatim}\ns

The most frequent structure of $\CH_{K_1}$ is
$\Z/8\Z \times \Z/8\Z \times\Z/2\Z \times\Z/2\Z$ with $272$ cases over $396$,
and there is never complete capitulation in $K_1/K$ (existence of elements of order $4$), 
but for $N$ large enough, capitulation (and even stability) is possible from $K_2$, 
as shown by the following examples with automatic capitulation in $K_3$:

\ft\begin{verbatim}
conductor f=20887 PK=x^3+x^2-6962*x-225889 CK=[4,4,2,2]
ell=193  N=5  Nn=2  n=1  CK=[4,4,2,2]  CK1=[8,8,2,2]
norm in K1/K of the component 1 of CK1: [2,4,0,0]~
norm in K1/K of the component 2 of CK1: [4,6,0,0]~
norm in K1/K of the component 3 of CK1: [4,0,0,0]~
norm in K1/K of the component 4 of CK1: [4,4,0,0]~

ell=193  N=5  Nn=2  n=2  CK=[4,4,2,2]  CK2=[8,8,2,2]
norm in K2/K of the component 1 of CK2: [4,0,0,0]~
norm in K2/K of the component 2 of CK2: [0,4,0,0]~
norm in K2/K of the component 3 of CK2: [0,0,0,0]~
norm in K2/K of the component 4 of CK2: [0,0,0,0]~
__________
ell=241  N=3  Nn=2  n=1  CK=[4,4,2,2]  CK1=[8,8,2,2]
norm in K1/K of the component 1 of CK1: [6,0,0,0]~
norm in K1/K of the component 2 of CK1: [0,6,0,0]~
norm in K1/K of the component 3 of CK1: [4,4,0,0]~
norm in K1/K of the component 4 of CK1: [4,0,0,0]~

ell=241  N=3  Nn=2  n=2  CK=[4,4,2,2]  CK2=[8,8,2,2]
norm in K2/K of the component 1 of CK2: [4,0,0,0]~
norm in K2/K of the component 2 of CK2: [0,4,0,0]~
norm in K2/K of the component 3 of CK2: [0,0,0,0]~
norm in K2/K of the component 4 of CK2: [0,0,0,0]~
__________
ell=1777  N=3  Nn=2  n=1  CK=[4,4,2,2]  CK1=[8,8,2,2]
norm in K1/K of the component 1 of CK1: [6,4,0,0]~
norm in K1/K of the component 2 of CK1: [4,2,0,0]~
norm in K1/K of the component 3 of CK1: [4,0,0,0]~
norm in K1/K of the component 4 of CK1: [0,4,0,0]~

ell=1777  N=3  Nn=2  n=2  CK=[4,4,2,2]  CK2=[56,56,2,2]
norm in K2/K of the component 1 of CK2: [4,0,0,0]~
norm in K2/K of the component 2 of CK2: [0,4,0,0]~
norm in K2/K of the component 3 of CK2: [0,0,0,0]~
norm in K2/K of the component 4 of CK2: [0,0,0,0]~
\end{verbatim}\ns

\smallskip
Then the following structures for $\CH_{K_1}$ are often obtained:

\smallskip
$\Z/4\Z \times \Z/4\Z \times \Z/4\Z \times \Z/4\Z$  ($42$ cases over $396$),

$\Z/4\Z \times \Z/4\Z \times \Z/2\Z \times \Z/2\Z \times \Z/2\Z \times \Z/2\Z$ 
($35$ cases over $396$),

$\Z/8\Z \times \Z/8\Z \times \Z/4\Z \times \Z/4\Z$  ($34$ cases over $396$).

\smallskip
Some cases of:

\smallskip
$\Z/8\Z \times \Z/8\Z \times\Z/2\Z \times\Z/2\Z \times\Z/2\Z \times \Z/2\Z$ ($8$ cases over $396$),

\smallskip\noindent
may give capitulation. Then we obtain a unique case of each of the following 
structures for $\CH_{K_1}$:

\smallskip
$\Z/4\Z \times \Z/4\Z \times \Z/4\Z \times \Z/4\Z \times \Z/2\Z \times \Z/2\Z$,

$\Z/8\Z \times \Z/8\Z \times \Z/4\Z \times \Z/4\Z \times \Z/2\Z \times \Z/2\Z$,

$\Z/8\Z \times \Z/8\Z \times \Z/8\Z \times \Z/8\Z$,

$\Z/16\Z \times \Z/16\Z \times \Z/2\Z \times \Z/2\Z \times \Z/2\Z \times \Z/2\Z$,

$\Z/32\Z \times \Z/32\Z \times \Z/2\Z \times \Z/2\Z \times \Z/2\Z \times \Z/2\Z$.

\smallskip
These examples suggest that for suitable $N$, capitulation is always
obtained in the tower in a wide variety of ways, and that each structure gives, 
most often, the same kind of results. These phenomena are certainly 
governed by precise probabilities.

\section{Numerical experiments over quadratic fields}
We consider real quadratic fields $K=\Q(\sqrt{-m})$ with $p \ne 2$. 
Let $L \subset K(\mu_\ell^{})$, $\ell \equiv 1 \pmod {2 p^N}$ inert in $K$, 
be a cyclic $p$-extension of degree $p^N$ of $K$. 
Of course, this case has no interest for verifications of the Main Conjecture
since it is true for the trivial reason $\chi = \varphi$; but it remains significant
to study the capitulation phenomenon. The program is analogous to the general 
one with slight modifications due to the quadratic context; we then give some 
excerpt for $p=3$ and $p=5$:

\subsection{General program for quadratic fields}\label{programquad}

\ft\begin{verbatim}
{p=3;N=2;Nn=2;bm=2;Bm=10^4;vHK=2;vHKn=2;Bell=500;
for(m=bm,Bm,if(core(m)!=m,next);PK=x^2-m;K=bnfinit(PK,1);
HK=K.no;if(valuation(HK,p)<vHK,next);CK=K.clgp;
forprime(ell=1,Bell,if(Mod(ell-1,2*p^N)!=0,next);
if(Mod(m,ell)==0,next);if(kronecker(m,ell)!=-1,next);
for(n=1,Nn,QKn=polsubcyclo(ell,p^n);P=polcompositum(PK,QKn)[1];
Kn=bnfinit(P,1);HKn=Kn.no;if(valuation(HKn,p)<vHKn,break);
htame=HKn/p^valuation(HKn,p);CKn=Kn.clgp[2];
print("PK=",PK," CK=",CK[2]," ell=",ell,
" N=",N," Nn=",Nn," n=",n," CK",n,"=",CKn);
rKn=matsize(CKn)[2];G=nfgaloisconj(Kn);Id=G[1];
for(k=2,2*p^n,Z=G[k];ks=1;
while(Z!=Id,Z=nfgaloisapply(Kn,G[k],Z);ks=ks+1);
if(ks==p^n,S=G[k];break));for(j=1,rKn,e=CKn[j];
A0=Kn.clgp[3][j];A=1;for(t=1,p^n,As=nfgaloisapply(Kn,S,A);
A=idealmul(Kn,A0,As));A=idealpow(Kn,A,htame);
X=bnfisprincipal(Kn,A)[1];
print("norm in K",n,"/K of the component ",j," of CK",n,":",X)))))}
\end{verbatim}\ns

\subsection{Case of quadratic fields and \texorpdfstring{$p=3$}{Lg}}

\subsubsection{Examples with \texorpdfstring{$n \in [1,2]$}{Lg}}  

\smallskip
For instance, the configuration:

\ft\begin{verbatim}
PK=x^2 - 1129 CK=[9]   
ell=13  N=1  Nn=1  n=1  CK=[9]  CK1=[27] 
norm in K1/K of the component 1 of CK1: [3]~
No capitulation in K1
\end{verbatim}\ns
since $\Nu_{K_1/K}(\CH_{K_1}) = \J_{K_1/K}(\CH_K)
= \CH_{K_1}^3$ of order $9$, there is no capitulation of $\CH_K$ in $K_1$. 

\ft\begin{verbatim}
PK=x^2 - 1129  CK=[9] 
ell=307  N=2  Nn=1  n=1  CK=[9]  CK1=[9,3] 
norm in K1/K of the component 1 of CK1:[3,0]~
norm in K1/K of the component 2 of CK1:[0,0]~
Incomplete capitulation in K1

ell=307  N=2  Nn=2  n=2  CK=[9]  CK2=[9,9]
norm in K2/K of the component 1 of CK2:[0,0]~
norm in K2/K of the component 2 of CK2:[0,0]~
Complete capitulation in K2
______
ell=19  N=2  Nn=1  n=1  CK=[9]  CK1=[9] 
norm in K1/K of the component 1 of CK1:[3]~
Incomplete capitulation in K1
Stability from K--->capitulation in K2

ell=19  N=2  Nn=2  n=2  CK=[9]  CK2=[9]
norm in K2/K of the component 1 of CK2:[0]~
Complete capitulation in K2
-----------------------------------------------
PK=x^2 - 1129  CK=[9]  
ell=73  N=2  Nn=2  n=1  CK=[9] CK1=[189,3] 
norm in K1/K of the component 1 of CK1: [84,0] 
norm in K1/K of the component 2 of CK1: [0,0] 
No capitulation in K1

ell=73  N=2  Nn=2  n=2  CK=[9]  CK2=[567,9] 
norm in K2/K of the component 1 of CK2: [252,0]
norm in K2/K of the component 2 of CK2: [0,0]  
No capitulation in K2
-----------------------------------------------
PK=x^2 - 3137 CK=[9] 
ell=199  N=2  Nn=2  n=1  CK=[9] CK1=[27,3] 
norm in K1/K of the component 1 of CK1: [21,0]~
norm in K1/K of the component 2 of CK1: [0,0]~  
No capitulation in K1

ell=199  N=2  Nn=2  n=2  CK=[9]  CK2=[27,9] 
norm in K2/K of the component 1 of CK2: [9,0]~
norm in K2/K of the component 2 of CK2: [0,0]~
Incomplete capitulation in K2
-----------------------------------------------
PK=x^2 - 8761 CK=[27]  
ell=19  N=2  Nn=2  n=1  CK=[27] CK1=[81] 
norm in K1/K of the component 1 of CK1: [3]~
No capitulation in K1

ell=19  N=2  Nn=2  n=2  CK=[27]  CK2=[81] 
norm in K2/K of the component 1 of CK2: [9]~
Incomplete capitulation in K2
\end{verbatim}\ns

We have the case of the following non-cyclic structure of $\CH_K$:

\ft\begin{verbatim}
PK=x^2 - 32009 CK=[3,3]  
ell=19  N=2  Nn=2  n=1  CK=[3,3]  CK1=[9,3] 
norm in K1/K of the component 1 of CK1:[3,0]
norm in K1/K of the component 2 of CK1:[0,0]
Incomplete capitulation in K1

ell=19  N=2  Nn=2  n=2   CK=[3,3]  CK2=[9,3]
norm in K2/K of the component 1 of CK2: [0,0]~
norm in K2/K of the component 2 of CK2: [0,0]~
Complete capitulation in K2
Stability from K1
-----------------------------------------------
PK=x^2 - 42817 CK=[3,3]
ell=19  N=2  Nn=2  n=1  CK=[3,3]  CK1=[9,3]
norm in K1/K of the component 1 of CK1: [3,0]~
norm in K1/K of the component 2 of CK1: [0,0]~

ell=19  N=2  Nn=2  n=2  CK=[3,3]  CK2=[27,3]
norm in K2/K of the component 1 of CK2: [9,0]~
norm in K2/K of the component 2 of CK2: [0,0]~
Incomplete capitulation in K2
\end{verbatim}\ns

\subsubsection{Examples of stability from \texorpdfstring{$K_1$}{Lg}}
Let's give the program testing the stability in $K_1/K$; as for the cubic 
case with $p=2$, there are much solutions. We put ${\sf ell = prime(t)}$,
${\sf t \in [2,nell]}$:

\ft\begin{verbatim}
{p=3;bm=2;Bm=10^4;N=3;nell=100;for(m=bm,Bm,if(core(m)!=m,next);
PK=x^2-m;K=bnfinit(PK,1);v=valuation(K.no,p);if(v!=N,next);
CK=K.clgp[2];for(t=1,nell,ell=prime(t);
if(t==nell,print("m=",m," nell insufficient");break);
if(Mod(ell-1,3)!=0,next);if(kronecker(m,ell)!=-1,next);
P=polcompositum(PK,polsubcyclo(ell,p))[1];K1=bnfinit(P,1);
v1=valuation(K1.no,p);if(v1!=N,next);
print("m=",m," PK=",PK," ell=",ell,
" CK=",K.clgp[2]," CK1=",K1.clgp[2]);break))}
N=2
m=1129 PK=x^2 - 1129 ell=19 CK=[9]  CK1=[9]
m=1654 PK=x^2 - 1654 ell=43 CK=[9]  CK1=[9]
m=3137 PK=x^2 - 3137 ell=19 CK=[9]  CK1=[9]
m=3719 PK=x^2 - 3719 ell=31 CK=[9]  CK1=[18,2]
(...)
m=9217 PK=x^2 - 9217 ell=7  CK=[18] CK1=[18]
m=9606 PK=x^2 - 9606 ell=31 CK=[18] CK1=[18,2,2]
m=9799 PK=x^2 - 9799 ell=19 CK=[18] CK1=[18,2,2]
m=9998 PK=x^2 - 9998 ell=37 CK=[9]  CK1=[9]
(...)
N=3
m=8761  PK=x^2 - 8761 ell=37  CK=[27] CK1=[27]
m=21433 PK=x^2 - 21433 ell=7  CK=[27] CK1=[27]
m=30859 PK=x^2 - 30859 ell=7  CK=[27] CK1=[27]
m=31327 PK=x^2 - 31327 ell=19 CK=[27] CK1=[27]
(...)
m=68513 PK=x^2 - 68513 ell=19 CK=[27] CK1=[27]
m=83713 PK=x^2 - 83713 ell=13 CK=[27] CK1=[54,2]
m=90271 PK=x^2 - 90271 ell=31 CK=[27] CK1=[54,2]
m=94865 PK=x^2 - 94865 ell=43 CK=[54] CK1=[54]
(...)
\end{verbatim}\ns

The most impressive is that, up to $m \leq 10^{10}$, small
primes $\ell$ are sufficient to get stability for cyclic groups $\CH_K$.

\subsection{Case of quadratic fields and \texorpdfstring{$p=5$}{Lg}}
We give some excerpt of numerical results, analogous to the case $p=3$;
most of examples give stability, whence capitulation in some layer.
We have taken ${\sf N=2, Nn=1}$.
For some rare cases, the capitulation is incomplete in the first layer.

\ft\begin{verbatim}
PK=x^2 - 24859 CK=[25]
ell=251  N=2  Nn=1  n=1  CK=[25]  CK1=[50,2,2,2]
norm in K1/K of the component 1 of CK1: [30,0,0,0]~
norm in K1/K of the component 2 of CK1: [0,0,0,0]~
norm in K1/K of the component 3 of CK1: [0,0,0,0]~
norm in K1/K of the component 4 of CK1: [0,0,0,0]~
Incomplete capitulation in K1
Stability from K--->capitulation in K_2

ell=401  N=2  Nn=1  n=1  CK=[25]  CK1=[1525]
norm in K1/K of the component 1 of CK1: [305]~
Incomplete capitulation in K1
-----------------------------------------------
PK=x^2 - 27689 CK=[25]
ell=101  N=2  Nn=1  n=1  CK=[25]  CK1=[25]
norm in K1/K of the component 1 of CK1: [5]~
Incomplete capitulation in K1
Stability from K--->capitulation in K_2
-----------------------------------------------
PK=x^2 - 68119 CK=[50]
ell=251  N=2  Nn=1  n=1  CK=[50]  CK1=[250]
norm in K1/K of the component 1 of CK1: [5]~
No capitulation in K1
-----------------------------------------------
PK=x^2 - 68819 CK=[25]
ell=101  N=2  Nn=1  n=1 CK=[25]  CK1=[25]
norm in K1/K of the component 1 of CK1: [5]~
Incomplete capitulation in K1
Stability from K--->capitulation in K_2

ell=151  N=2  Nn=1  n=1  CK=[25]  CK1=[125]
norm in K1/K of the component 1 of CK1: [5]~
No capitulation in K1
-----------------------------------------------
PK=x^2 - 69403 CK=[25]
ell=251  N=2  Nn=1  n=1  CK=[25]  CK1=[25,5]
norm in K1/K of the component 1 of CK1: [5,0]~
norm in K1/K of the component 2 of CK1: [0,0]~
Incomplete capitulation in K1

ell=401  N=2  Nn=1  n=1  CK=[25]  CK1=[50,2,2,2]
norm in K1/K of the component 1 of CK1: [30,0,0,0]~
norm in K1/K of the component 2 of CK1: [0,0,0,0]~
norm in K1/K of the component 3 of CK1: [0,0,0,0]~
norm in K1/K of the component 4 of CK1: [0,0,0,0]~
Incomplete capitulation in K1
Stability from K--->capitulation in K_2
-----------------------------------------------
PK=x^2 - 88211 CK=[25]
ell=101  N=2  Nn=1  n=1  CK=[25]  CK1=[125]
norm in K1/K of the component 1 of CK1: [5]~
No capitulation in K1

ell=151  N=2  Nn=1  n=1  CK=[25]  CK1=[25]
norm in K1/K of the component 1 of CK1: [5]~
Incomplete capitulation in K1
Stability from K--->capitulation in K_2
-----------------------------------------------
PK=x^2 - 119029 CK=[50]
ell=251  N=2  Nn=1  n=1  CK=[50]  CK1=[250]
norm in K1/K of the component 1 of CK1: [5]~
No capitulation in K1

ell=401  N=2  Nn=1  n=1  CK=[50]  CK1=[50]
norm in K1/K of the component 1 of CK1: [5]~
Incomplete capitulation in K1
Stability from K--->capitulation in K_2
\end{verbatim}\ns

\section{On the use of the auxiliary extensions \texorpdfstring{$K(\mu_\ell^{})/K$}{Lg}}

\subsection{Remarks on the classical proof of the real Main Conjecture}
Contrary, to our point of view, some proofs of the Main Conjecture, 
as that of Thaine \cite{Thai1988}, have a fundamental difference.
The brief overview that we give now must be completed by technical elements 
that the reader can find especially in \cite[\S\S\,15.2, 15.3]{Was1997}.

\smallskip
Let $f$ be the conductor of $K$. Let $\ell \equiv 1 \pmod {2p^N}$, 
{\it totally split in $K$} and $p^N \geq p^e$, the exponent of $\CH_K$; 
put $M_0 = \Q(\mu_\ell^{})$ and $M := M_0 K$. From Proposition \ref{cycloformula}, 
the cyclotomic unit $\eta_M^{}$ of $M$ is such that $\Norm_{M/K}(\eta_M^{}) = 1$
and the link with $\eta_K^{}$ is only a congruence modulo $(1-\zeta_\ell)$.

\smallskip
Put $\eta_M^{} = \alpha_M^{1-\sigma}$ where
$\alpha_M^{} \in M^\times$ is such that 
$(\alpha_M^{}) \in I_M^G$; modulo $K^\times$, we can take
$\alpha_M^{}$ integer in $M$, whence:
\begin{equation}\label{resolvent}
(\alpha_M^{}) = \jj_{M/K}^{}({\mathfrak a}_K) \cdot {\mathfrak L}_1^{\Omega_\ell}, 
\end{equation}  
where ${\mathfrak a}_K \in I_K$ and ${\mathfrak L}_1$ is a fixed prime ideal dividing 
$\ell$ in $M$, thus totally ramified in $M/K$, with
${\Omega_\ell} = \sm_{s \in g} r_s \cdot s^{-1}$, $r_s \geq 0$;
thus, since $\Norm_{M/K}({\mathfrak L}_1) = {\mathfrak l}_1$, 
${\mathfrak L}_1 \mid {\mathfrak l}_1 \mid \ell$ in $M/K/\Q$, this yields
$(\alpha_K) := (\Norm_{M/K}(\alpha_M^{})) = 
{\mathfrak a}_K^{\ell-1} \cdot {\mathfrak l}_1^{\,{\Omega_\ell}}$.
But ${\mathfrak a}_K^{\ell - 1}$ is principal since $\ell -1 \equiv 0 \pmod {p^e}$, 
whence ${\mathfrak l}_1^{\Omega_\ell}$ principal.

\begin{lemma}\label{nontrivial}
Except a finite number of primes $\ell$, the ideal ${\mathfrak L}_1^{\Omega_\ell}$ 
in relation \eqref{resolvent} gives a non trivial relation, in the meaning that 
${\Omega_\ell}$ is not of the form $\lambda \cdot \Nu_{M/M_0}$, 
$\lambda \geq 0$, giving ${\mathfrak l}_1^{\Omega_\ell} = (\ell)^\lambda$.
\end{lemma}

\noindent{\bf Proof.}
Assume that ${\Omega_\ell} = \lambda \cdot \Nu_{M/M_0}$; the character of
${\mathfrak L}^{\Omega_\ell}_1$, as $\Z[g]$-module,
is the unit one and any non-trivial 
$\varphi$-component $\alpha_{M,\varphi}^{}$ of $\alpha_{M}^{}$ is
prime to $\ell$, thus congruent, modulo any ${\mathfrak L} \mid \ell$ in $M$, to 
$\rho_{\mathfrak l}^{} \in \Z$, $\rho_{\mathfrak l}^{} \not\equiv 0 \pmod \ell$ 
(residue degrees~$1$ in $M/\Q$). 
Since ${\mathfrak L}^\sigma = {\mathfrak L}$, we obtain $\eta_{M,\varphi^{}}^{} 
=\alpha_{M,\varphi}^{\sigma-1} \equiv 1 \pmod {\mathfrak L}$.

\smallskip
We have $\eta_{f\ell}^{} \equiv \eta_f^{} \pmod {(1-\zeta_\ell)}$ where
$1-\zeta_\ell$ is an uniformizing parameter at the places above $\ell$ in 
$M_0$, so that $\eta_M^{} \equiv \eta_K^{} \pmod{(1-\zeta_\ell)}$, which leads
to $\eta_{K,\varphi^{}}^{} \equiv \eta_{M,\varphi} \equiv 1 \pmod {\mathfrak l}$, for all 
${\mathfrak l} \mid \ell$, giving $\eta_{K,\varphi^{}}^{} \equiv 1 \pmod \ell$ 
(absurd for almost all $\ell$).
\qed

\medskip
Reducing modulo $\Nu_{M/M_0}$, one may get ${\Omega_\ell} \ne 0$,
``minimal'' in an obvious sense, with $r_s \geq 0$ but not all zero.

\smallskip
Consider $\ds \frac{\alpha_M^s}{(1-\zeta_\ell)^{r_s}}$ modulo ${\mathfrak L}_1$ 
and the conjugations $\alpha_M^\sigma = \alpha^{}_M \cdot \eta^{}_M$
and $\ds \frac{1-\zeta_\ell^\sigma}{1-\zeta_\ell} = 
\frac{1-\zeta_\ell^{\g_\ell^{}}}{1-\zeta_\ell} \equiv \g_\ell^{} \pmod{(1-\zeta_\ell)}$ 
(where $\g_\ell^{}$ is a primitive root modulo $\ell$ such that $\zeta_\ell^\sigma =: 
\zeta_\ell^{\g_\ell^{}}$); one gets:
\[\Big( \frac{\alpha_M^s} {(1-\zeta_\ell)^{r_s}} \Big)^\sigma \!\!=
 \frac{\alpha_M^{\sigma s}}{(1-\zeta_\ell)^{\sigma r_s}} 
\equiv \frac{\eta_M^s \alpha_M^s}{(\g_\ell (1-\zeta_\ell))^{r_s}} \equiv
\frac{\eta_M^s} {\g_\ell^{r_s}}  
\Big( \frac{\alpha_M^s} {(1-\zeta_\ell)^{r_s}} \Big)\!\!\!\pmod{{\mathfrak L}_1}, \]
whence $\g_\ell^{r_s} \equiv \eta_M^s \equiv \eta_K^s \pmod {{\mathfrak l}_1}$, 
which identifies the coefficients $r_s$.

\smallskip
So, this yields a non-trivial relation between the classes of the 
conjugates of ${\mathfrak l}_1$ that is $p$-localizable; this constitute the basis 
of the reasonings, on condition to add many more technical arguments to get some
annihilation of $\CE_{K,\varphi}/\CF_{K,\varphi}$, then a final equality 
$\order \CH_{K,\varphi} = (\CE_{K,\varphi} : \CF_{K,\varphi})$.

\smallskip
In this way, Thaine's method is essentially analytic, working on norm properties and 
subtle congruences of the cyclotomic units, leading to the principle of Kolyvagin Euler 
systems, while that using capitulation (if any) is of class field theory framework,
especially with Chevalley--Herbrand context,
and gives immediately the result {\it without any supplementary work}.

\subsection{Heuristics about the Conjecture \ref{conjcap} of capitulation}
Suppose, in the totally opposite case, that for all inert primes $\ell \equiv 1\!\pmod {2 p^N}$, 
with $N$ arbitrary large, there is never any capitulation in the real cyclic $p$-tower $L/K$ of 
$K(\mu_\ell^{})/K$; we get easily the following result, where we denote by $K_n$ 
the subfield of $L$ of degree $p^n$ over $K$, $n \in [0, N]$.

\begin{proposition}
Assume that, for all $n, m \in [0,N]$, $n \leq m$,
the transfer maps $\J_n^m := \J_{K_m/K_n}$ are injective.
Put $G_n^m := \Gal(K_m/K_n)$, $\CH_n := \CH_{K_n}$, $\CH_m := \CH_{K_m}$.
Then $\order \CH_m \geq \order \CH_n \cdot \order \CH_n[p^{m-n}]$ which 
leads, for all $n$, to $\order \CH_{K_n} \geq \order \CH_K \cdot p^{n \cdot r_K}$
where $r_K^{} = \dim_{\F_p}(\CH_K/\CH_K^p)$ is the $p$-rank of $\CH_K$.
\end{proposition}

\noindent{\bf Proof.}
From the exact sequence of $\Z[G_n^m]$-modules:
\[1 \to  \J_n^m\CH_n \to \CH_m \to \CH_m/ \J_n^m \CH_n \to 1, \] 
we get, from the Chevalley--Herbrand formula 
$\order \CH_m^{G_n^m} = \order \CH_n$ and the injectivity of $\J_n^m$,
$\CH_m^{G_n^m} / \J_n^m \CH_n = 1$ and the exact sequence:
$$1  \to (\CH_m/ \J_n^m \CH_n)^{G_n^m}  
 \to  {\rm H}^1(G_n^m, \J_n^m \CH_n) \to  {\rm H}^1(G_n^m, \CH_m), $$

\noindent
where
${\rm H}^1(G_n^m, \J_n^m \CH_n) = (\J_n^m \CH_n) [p^{m-n}] \simeq \CH_n[p^{m-n}]$
and:
$$\hspace{-5.7cm}\order {\rm H}^1(G_n^m, \CH_m) = \order {\rm H}^2(G_n^m, \CH_m)$$
$$\hspace{2.0cm}= \order \CH_m^{G_n^m} /\order \J_n^m \circ \Norm_n^m(\CH_m) =
\order \CH_m^{G_n^m} /\order \CH_n = 1$$ 
giving $(\CH_m/ \J_n^m \CH_n)^{G_n^m} \simeq \CH_n[p^{m-n}]$. 
Whence $\order \CH_m \geq \order \J_n^m \CH_n \cdot \order \CH_n[p^{m-n}]$ $ = 
\order \CH_n  \cdot \order \CH_n[p^{m-n}]$.
For $m=n+1$, one obtains $\order \CH_{n+1}\geq \order \CH_n \cdot p^{\rk_p(\CH_n)}$,
then the last claim by induction, ${\rm rank}_p(\CH_n)$ being increasing.
\qed

\medskip
This result indicates that, for $L=K_N$, the filtration $(\CH_L^i)_{i \geq 0}$
has length unbounded regarding $N$. Since $\order (\CH_L^{i+1}/\CH_L^i) =
\ds \frac{\order \CH_K}{\order \Norm_{L/K}(\CH_L^i)}$, giving a decreasing sequence,
in an probabilistic point of view, one sees that the length of the filtration depends 
essentially on the size of $\CH_K$ but not necessarily of~$N$.

\subsection{Conclusion}
The behavior of $p$-class groups in real cyclic $p$-towers $L/K$
(with $L \subset K(\mu_\ell^{})$, $\ell \equiv 1 \pmod {2 p^N}$
inert in $K/\Q$), suggests that there exist infinitely many 
such primes for which $\CH_K$ capitulates at some layer.
The most interesting fact being that, if so, this implies the real Main Conjecture 
on abelian fields in the semi-simple case with minimal analytic arguments
and almost trivial proof from the classical context of Chevalley--Herbrand 
formulas.

\smallskip
Regarding some works about these questions of ``exceptional classes'' 
(i.e., non invariant), it is admitted (and proved in some circumstances)
that the filtrations defining the $\CH_{K_n}$, for $n \in [1, N]$, are random
and that cases of ``unbounded'' filtrations are of probabilities tending to $0$ with $N$.
For instance, in \cite{KoPa2022,Smi2022} it is proved that $p$-class groups of cyclic 
extensions $L/\Q$ of degree $p$ have a standard distribution, the case of filtrations
of length $1$ (i.e., $\CH_L = \CH_L^G$) being the most probable; but this is another
context since the Chevalley--Herbrand formula $\order \CH_L^G = p^{r-1}$ is non-trivial
only if the number $r$ of ramified primes is at least $2$, in which case
$\CH_L^G = \CH_L^\ram$ (see Theorem \eqref{suitecap}) and $\order( (\CH_L/\CH_L^G)^G) 
=  \ffrac{p^{r-1}}{(\Lbda : \Lbda \cap \Norm_{L/\Q}(L^\times))}$ (see \eqref{filtration}), 
where $\Lbda$ is generated by the $r$ prime numbers ramified in $L/\Q$; 
whence the algorithm giving the $\order(\CH_L^{i+1}/\CH_L^i) =  
\ffrac{p^{r-1}}{(\Lbda_i : \Lbda_i \cap \Norm_{L/\Q}(L^\times))}$.
But, in our case ($\CH_K \ne 1$ while $r=1$), on the contrary, the norm factors are 
trivial and $\order(\CH_L^{i+1}/\CH_L^i) = \ffrac{\order \CH_K}{\Norm_{L/K}(\CH_L^i)}$
depends of the structure of $\CH_K$ and allows standard probabilities about order and
structure of $\CH_L$.

\smallskip
In other words, the context of ``fixed points formulas'' in cyclic $p$-towers, is, in some 
sense, specific and ``easier'' or ``more reachable'' than the general framework on 
$p$-class groups in arbitrary number fields, only accessible via complex analytic
methods, as the reader can see in the recent papers \cite{EPW2017, Wan2020, 
PTBW2020, Gra2020, KlWa2022, Gra2022$^b$, Pier2022}, among many others. 

\smallskip
Moreover, the $\epsilon$-conjectures are of no help, here, since the discriminants of the fields
$K(\mu_\ell^{})$ are larger than $\ell^{\ell - 2}  = O(p^{N \cdot p^N})$ giving
dreadful bounds for the $\order \CH_L $, and the philosophy is on the contrary (as for 
Greenberg's conjecture in $L=K_\infty$ for ``$\ell = p$'') that $\order \CH_L = 
O(\order \CH_K)$ for infinitely many $\ell$'s.
To be more audacious, one may imagine an ``Iwasawa behavior'' leading, for $N \gg 0$, 
to formulas of the form $\order \CH_{K_n} = p^{\lambda(\ell) \cdot n + \nu(\ell)}$ for 
$n \in [n_0(\ell), N]$, possibly with $\lambda(\ell) = 0$ (stability) for infinitely many $\ell$'s.

\smallskip
In other words, capitulation phenomena in $p$-towers are of $p$-adic type in a 
specific meaning adding randomness and, logically, govern many arithmetic results 
and conjectures of number theory around class field theory, and appear essentially 
as new general regularity principle which deserves deepened researches.

\end{document}